\documentclass[11pt]{article}
\usepackage[english]{babel}
\usepackage[cp850]{inputenc}
\usepackage[dvips]{graphicx}
\usepackage{amsmath,amsfonts,amsthm,amssymb}
\usepackage[usenames, dvipsnames]{color}
\usepackage{fancyhdr}
\usepackage{stmaryrd}
\usepackage[colorlinks=true,citecolor=red,linkcolor=blue,urlcolor=blue,pdfstartview=FitH]{hyperref}
\usepackage{dsfont}
\usepackage{xcolor}
\usepackage{epsfig}

\font\tenmath=msbm10 scaled 1200

\font\sevenmath=msbm7 scaled 1200
\font\fivemath=msbm5 scaled 1200

\newfam\mathfam \textfont\mathfam=\tenmath
\scriptfont\mathfam=\sevenmath \scriptscriptfont\mathfam=\fivemath

\def\R{{\mathbb R}}
\def\N{{\mathbb N}}
\def\E{{\mathbb E}}

\def\P{{\mathbb P}}

\newtheorem{Thm}{Theorem}[section]
\newtheorem{Lem}{Lemma}[section]
\newtheorem{Pro}{Proposition}[section]
\newtheorem{Cor}{Corollary}[section]

\newfam\mathfam \textfont\mathfam=\tenmath
\scriptfont\mathfam=\sevenmath \scriptscriptfont\mathfam=\fivemath

\def \^#1{\if#1i{\accent"5E\i}\else{\accent"5E#1}\fi}

\def \a{\alpha}

\def \n{\eta}

\def \ds{\displaystyle}
\def \cqfd{\quad\Box}

\def \ss{\smallskip}
\def \bs{\bigskip}
\def \ni{\noindent}

\def\ni{\noindent}

\title{\bf Greedy vector quantization}

\author{
\textsc{Harald Luschgy}~\thanks{Universit\"at Trier, FB IV-Mathematik, D-54286 Trier, Germany. E-mail: {\tt luschgy@uni-trier.de}}  $\quad$ and $\quad$ 
\textsc{Gilles Pag\`es} \thanks{Laboratoire de Probabilit\'es et Mod\`eles Al\'eatoires, UMR~7599, UPMC, case 188, 4, pl. Jussieu, F-75252 Paris Cedex 5, France. E-mail: \texttt{gilles.pages@upmc.fr}} }

\date{}

\begin{document}

\maketitle

\begin{abstract} We investigate  the greedy version of the $L^p$-optimal vector quantization problem for an  $\R^d$-valued random vector $X\!\in L^p$. We show the existence of a sequence $(a_N)_{N\ge 1}$ such that $a_N$ minimizes $a\mapsto\big \|\min_{1\le i\le N-1}|X-a_i|\wedge |X-a|\big\|_{L^p}$ ($L^p$-mean quantization error at level $N$ induced by $(a_1,\ldots,a_{N-1},a)$).  We show that this sequence produces $L^p$-rate optimal $N$-tuples $a^{(N)}=(a_1,\ldots,a_{_N})$   ($i.e.$ the $L^p$-mean quantization error  at level $N$ induced by $a^{(N)}$  goes to $0$ at rate $N^{-\frac 1d}$). Greedy optimal sequences also  satisfy, under natural additional assumptions,  the distortion mismatch property: the $N$-tuples $a^{(N)}$ remain rate optimal with respect to the $L^q$-norms, $p\le q <p+d$. Finally, we  propose optimization methods to compute greedy sequences, adapted from usual Lloyd's~I and Competitive Learning Vector Quantization procedures, either in their deterministic (implementable when $d=1$)  or  stochastic versions.
\end{abstract}


\bigskip

\ni {\em 2010 AMS Classification:}  60G15, 60G35, 41A25.

\section{Introduction and definition  of greedy quantization sequences}

Let  $p\!\in (0,+\infty)$ and $L_{\R^d}^p(\Omega,{\cal A},\P)=\{Y:(\Omega,{\cal A},\P)\to \R^d$,\; measurable, $\|Y\|_p = \big(\E|Y|^p\big)^{\frac 1p}<+\infty\}$ where $|\,.\,|$ denotes a norm  on $\R^d$. We consider $X:(\Omega,{\cal A},\P)\to \R^d$ an $L^p$-integrable random vector. For every  $\Gamma\subset \R^d$, we define the $L^p$-mean quantization error induced by $\Gamma$ as the  $L^p$-mean of the distance of the random vector $X$ to the subset $\Gamma$ (with respect to  the norm $|\,.\,|$), namely
\[
e_p(\Gamma,X)= \big\|d(X,\Gamma)\big\|_p
\] 
where $d(\xi,A)= \inf_{a\in A} |\xi-a|$, $\xi\!\in \R^d$, $A\subset \R^d$, denotes the distance of $\xi $ to $A$. This quantity is always finite when $X\!\in L^p(\P)$  since $e_p(\Gamma,X)\le \|X\|_p+ \min_{a\in \Gamma}|a|<+\infty$ owing to Minkowski's inequality when $p\ge 1$.  When $p\!\in (0, 1)$, one has likewise $e_p(\Gamma,X)^p\le \|X\|^p_p+ \min_{a\in \Gamma}|a|^p<+\infty$. The usual $L^p$-optimal quantization problem {\em at level} $ N\ge 1$  is to solve the following minimization problem 
\begin{equation}\label{eq:epNX}
e_{p,N}(X)= \min_{\Gamma\subset \R^d, |\Gamma| \le N} e_p(\Gamma,X)
\end{equation}
where $|\Gamma|$ denotes the cardinality of the subset $\Gamma $, sometimes called {\em grid} in Numerical Probability or {\em codebook} in Signal processing. The use of ``$\min$" instead of ``$\inf$" is justified by the fact (see~Proposition~4.12 in~\cite{GRLU}, p.47 or~\cite{PAG0}) that  this infimum is always attained by an {\em optimal quantization} grid $\Gamma^{(N)}$ (of full size $N$ if the support of the distribution $\mu=\P_{_X}$ of $X$ has at least $N$ elements).

The above optimal vector quantization problem is clearly related to the approximation rate  of an $\R^d$-valued random vectors $X:(\Omega,{\cal A}, \P)\to \R^d$ by random vectors taking at most $N$ values ($N\!\in \N)$. One shows  (see $e.g.$ Theorem~4.12 in~\cite{GRLU} combined with comments, Section~3.3, p.33) that, for very $p\!\in (0,+\infty)$,
\begin{eqnarray*}
e_{p,N}(X) &=& \min\Big\{\|X- q(X)\|_{_p}, \,q:\R^d\to \R^d, \mbox{ Borel}, \,|q(\R^d)|\le N\Big\}\\
&=& \min\Big\{\|X- Y\|_{_p}, \,Y:\Omega \to \R^d, \mbox{ measurable}, \,|Y(\Omega)|\le N\Big\},
\end{eqnarray*}
both minima being attained by random vectors of the form
\begin{equation}\label{eq:minVor}
Y^{(N)}=\widehat X^{(N)}:=\pi_{\Gamma^{(N)}}(X)
\end{equation}
where $\pi_{\Gamma^{(N)}}$ denotes a {\em Borel projection on $\Gamma^{(N)}$ following the nearest neighbour rule} where $\Gamma^{(N)}\subset \R^d$ has size at most $N$.

This modulus is also related to the Wasserstein (pseudo-)distance ${\cal W}_p$, $p\!\in (0,1]$ on the space of Borel probability measure on $\R^d$: let ${\cal P}_{_N}$ be the set of distributions whose support has at most $N$ elements.   Let $\mu$ be a Borel distribution on $\R^d$ and let $ \nu\!\in {\cal P}_{_N}$ that we can associate to   random vectors $X$ and $Y$ respectively ; then for every $p$-H\" older function $f:\R^d\to \R$, with $p$-H\"older ratio $[f]_{_{p, Hol}}<+\infty$ and every $ \nu\!\in {\cal P}_{_N}$, 
\begin{equation}\label{eq:cubat}
\Big|\int_{\R^d}f\, d\mu- \int_{\R^d}f\, d \nu\Big| = \Big|\E\, f(X)-\E\, f(Y)\Big|\le [f]_{_{p, Hol}}\|X-Y\|_{_p}.
\end{equation}
Conversely, noting that the function $\xi \mapsto d(\xi, \Gamma^{(N)})$ is $p$-H\"older, we easily derive that 
\[
{\cal W}_p(\mu,{\cal P}_{_N}) =\inf_{\nu \in {\cal P}_{_N}} \sup\Big\{\Big | \int_{\R^d}f\, d\mu- \int_{\R^d}f\, d \nu\Big|,\;[f]_{_{p, Hol}}\le 1\Big\}= e_{p,N}(X)
\]
 When $\nu= \mu \circ  \pi^{-1}_{\Gamma^{(N)}}={\cal L}(Y^{(N)})$ (defined in~\eqref{eq:minVor}), the above inequality~\eqref{eq:cubat} is often used as a cubature formula for numerical integration (see~\cite{PAG0, Cher1, Cher2}). When dealing directly with  with random vectors, extensions of this formula are used to compute  conditional expectations (see~among others~\cite{BAPA2, PHSERU, BRDeSAPDU} and further on for more references).

\smallskip
The most celebrated result in Optimal (Vector) Quantization Theory is undoubtedly Zador's Theorem  (see~\cite{ZA,BUWI} and \cite{GRLU}) recalled below which rules the {\em sharp} asymptotic rate of convergence of $e_{p,N}(X)$ as the {\em quantization level} $N$ (or grid size)  goes to infinity.

\begin{Thm}[(Zador's Theorem), see~\cite{GRLU}, Theorem~6.2, p.78 and Remark~6.3$(c)$, see also~\cite{GrLuPa2}] \label{thm:Zador}$(a)$ If  $\E|X|^{p}<+\infty$ and $\mu= \P_{X}=  \varphi. \lambda_d + \nu$ where $\nu $ is a singular Borel measure with respect to the Lebesgue measure $\lambda_d$ on $\R^d$.  Then
\[
\liminf_N N^{\frac 1d}  e_{p,N}(X) \ge \widetilde J_{p,d}\left(\int_{\R^d} \varphi^{\frac{d}{p+d}}d\lambda_d\right)^{\frac 1p+\frac 1d}
\]
where $\widetilde J_{p,d}$ is the sharp limit for the uniform distribution $U([0,1]^d)$ over the unit hypercube which  satisfies 
\[
\widetilde J_{p,d} = \inf_N N^{\frac 1d}  e_{p,N}\big(U([0,1]^d)\big) \!\in (0,+\infty).
\]
\noindent $(b)$ If furthermore  $\E|X|^{p+\delta}<+\infty$ or some $\delta>0$, then 
\[
\lim_N N^{\frac 1d}  e_{p,N}(X) = \widetilde J_{p,d}\left(\int_{\R^d} \varphi^{\frac{d}{p+d}}d\lambda_d\right)^{\frac 1p+\frac 1d}.
\]
\end{Thm}

This $N^{-\frac 1d}$ (sharp) rate is known as {\em the curse of dimensionality}. The numerical search of optimal grids solution to~\eqref{eq:epNX} (especially in the quadratic setting when $d=2$) leads to an $N\times d$-dimensional problem for each grid size $N$ which is often too demanding in practice when $N$ or $d$ grows. Hence the need for a possibly sub-optimal ``solution"" to this problem, easier to compute in terms of complexity and dimensionality, provided the price to pay remains asymptotically reasonable. 

\bigskip
The starting idea of {\em greedy quantization} is to determine a {\em sequence } $(a_N)_{N\ge 1}$ of points of $\R^d$ which is  recursively optimal {\em step by step} or {\em level by level} with respect to the $L^p$-mean quantization criterion. We mean that, if we set $a^{(N)} = \{a_1,\ldots,a_{_N}\}$, $N\ge 1$, and $a^{(0)}=\emptyset$, then 
\begin{equation}\label{eq:greedy}
\forall\, N\ge 0, \quad a_{N+1}\!\in {\rm argmin}_{\xi\in \R^d} e_p(a^{(N)}\cup\{\xi\}, X).
\end{equation}
 Note that  $a_1$ is simply an  $L^p$-median of (the distribution of) $X$ and that, when $p>1$,  a strict convexity argument implies the uniqueness of this $L^p$-median.  This  idea to design not only optimal $N$-tuples  but an optimal sequence which, hopefully, will produce   $N$-tuples with a rate optimal behavior as $N\to+\infty$  is very natural and can be compared to  sequences  with low discrepancy in Quasi-Monte Carlo methods. 
 
 In fact, such sequences have already been investigated  in an $L^1$ setting for compactly supported random vectors $X$  as a model of short term experiment planning $vs$ long term experiment planning represented by regular optimal quantization at a given level $N$ (see~\cite{Brancoetal}). Our aim in this paper is to solve this {\em greedy}  optimization problem for as general  as possible distributions $\mu=\P_{_X}$ and in any $L^p$-space, $p\!\in(0,+\infty)$, in two directions: first establish the existence of such {\em $L^p$-optimal greedy sequences} and then evaluate their  rate of decay of $e_p(a^{(N)},X)$ to $0$ as the quantization level $N$ goes to infinity.
 
 A possible wider field of applications is to substitute such sequences to optimal $N$-quantizers in the quantization based  numerical	schemes that have been developed in the early 2000s. In these procedures optimal quantizations used as a spatial discretization method that ``fits" optimal the distribution of interest at each time step. Among these application, often in connection with  Finance but also with  reliability, we may mention Numerical integration (see~\cite{PAG0,PAPR}), Optimal Stopping Theory (pricing of American style or callable derivatives, see~\cite{BAPA1,BAPA2,BAPA3}), Stochastic control of diffusions and portfolio optimization (see~\cite{PAPH, PHSERU, COPHRU}), or control of PDMP(\footnote{Piecewise Constant Deterministic Markov Processes introduced by M. Davis in~\cite{DAV}.}),  for reliability (see~\cite{BRDeSAPDU, DeSAPDU}), non-linear filtering and stochastic volatility models~(see~\cite{PAPH}), discretization of $BSDE$s and Stochastic $PDE$s~'see~\cite{GOPAPHPR}). See also the review papers~\cite{PAPHPR,PAPR3} and the references therein for more details.
 In most of  these applications, up to some variant, an $\R^d$-valued discrete time Markov chain $(X_k)_{0\le k\le n}$  is approximated path wise and    in distribution  by its quantized approximation sequence $(\widehat X^{\Gamma_k}_k)_{0\le k\le n}$   living on a {\em quantization tree} made up by the  {\em optimal quantization grids} $\Gamma_k$ (of varying  sizes $N_k$) and  the transitions matrices $\pi^k={\cal L}\big(\widehat X^{\Gamma_{k+1}}\,|\, \widehat X^{\Gamma_k}\big)$ which discretize the Markov dynamics   of the chain. The quantization based scheme turns out to be in many cases  spatial discretization of a (Backward) Dynamical Programming principle.
 Given the common sizes  of the  grids in these  implemented procedures ($N_k$ is often greater than $1\, 000$) and the number $n$ of time steps  ($n\ge 10$ and sometimes equal to $100$) the storing of this quantization tree may exceed the storage capacity of the computing device. Using  the induced grids $a^{(N_0)},\, a^{(N_1)}\dots,a^{(N_n)}$ induced by a greedy optimal sequence $(a_N)_{N\ge 1}$  will dramatically reduce this drawback, provided that, on the other hand, their  rate of decay of their mean quantization rates remain comparable to those of optimal quantizers.
 
 \medskip The paper is organized as follows: in Section~\ref{sec:Exist}, the existence of $(L^p,\mu)$-optimal greedy  sequences and their first properties are established for general and  Euclidean norms. In Section~\ref{sec:Rate}, $(L^p,\mu)$-optimal greedy  sequences are shown to be rate optimal in terms of  mean quantization error, compared to sequences of $L^p(\mu)$-optimal $N$-quantizers. We also solve~-~positively~--~the so-called {\em distortion mismatch problem} $i.e.$ the property that the above rate optimal decay property remains true for  the $L^q(\mu)$-mean quantization error when $q\!\in [p, p+d)$ in a $d$-dimensional setting (and sometimes for $q=p+d$).   In Section~\ref{sec:PractiCriterion}, easy-to-check criteria, mostly borrowed   from~\cite{GrLuPa2}, are adapted to our greedy framework. Section~\ref{sec:Few} is devoted to some further questions about the asymptotic behaviour of $L^p$-greedy sequences, compared to $L^p$-optimal $N$-quantizers or non-greedy $L^p$-rate optimal sequences.  In Section~\ref{sec:greedyalgo}, we propose numerical procedures to compute quadratic optimal greedy  sequences  in both $1$ and higher dimensional settings,  either by deterministic means or by simulation. Finally, we propose in Section~\ref{sec:GreedyvsQMC}, when $X$ is uniformly distributed on the unit hypercube $[0,1]^d$, a  comparison between optimal greedy sequences and  the   {\em sequences with low discrepancy}  popularized by the Quasi-Monte Carlo method.


\medskip
\noindent {\sc Notations:}  $\bullet$   $\N^*=\{1,2,\ldots\}$ the set of positive integers.

\smallskip
\noindent $\bullet$ $|\,.\,|$ denotes any norm on $\R^d$ (except specific mention). $B(x,\rho)$ denotes the closed ball centered at $x\!\in \R^d$ with radius $\rho>0$. For  every subset $A\subset \R^d$ and $\xi\!\in \R^d$,  $d(\xi, A) = \inf_{a\in A} |\xi-a|$ (distance of $\xi$ to the set $A$ in $(\R^d,|\,.\,|)$).  

\section{Existence of optimal greedy  quantization sequences}\label{sec:Exist}

\begin{Pro}\label{pro:greedex} $(a)$ {\em Existence:}   If $X\!\in L^p(\P)$, then the sequence of optimization problems~\eqref{eq:greedy} admits at least one solution $(a_{_N})_{N\ge 1}$ where $a_1$ is the $L^p$-median of the distribution $\mu$. Moreover, the finite sequence $\big(e_p(a^{(n)},X)\big)_{1\le n\le N}$ is (strictly) decreasing  as long as $N\le {\rm card}\big({\rm supp}(\mu)\big)$. In particular, $a_n\notin a^{(n-1)}$, for $n\!\in\{1,\ldots,N\}$. 

Any such a solution is called an $L^p$-{\em optimal  greedy quantization  sequence}. 

\smallskip
\noindent $(b)$ {\em Local optimality:} As long as $N\le {\rm card}\big({\rm supp}(\mu)\big)$ 
$$
\mu\big(C_{a_N}(a^{(N)})\big)>0\,\mbox{ where} \; C_{a_N}(a^{(N)})=\big\{\xi\!\in \R^d\,|\, |\xi-a_{_N}|< \min_{1\le i\le N-1}|\xi-a_i|\big\}
$$
and for any Borel set $C$ such that $C_{a_N}(a^{(N)})\subset C \subset W_{a_N}(a^{(N)})=\big\{\xi\!\in \R^d\,|\, |\xi-a_{_N}|\le  \min_{1\le i\le N-1}|\xi-a_i|\big\}$, 
$a_{_N}$ is solution to the local optimization problem 
\[
a_{_N} \in {\rm argmin}_{a\in \R^d} \E\big(|X-a|^p\,|\, X\in C\big).
\]

\smallskip
\noindent $(c)$ {\em Space filling:} Assume $X\!\in L^q_{\R^d}(\P)$ for some $q\ge p$. Then, any  $L^p$-optimal  greedy quantization  sequence $(a_N)_{N\ge 1}$   satisfies  
\[
 \lim_N e_q(a^{(N)},X)=0
\] 
$i.e.$, equivalently,   $\displaystyle \lim_{N\to+\infty } \int_{\R^d} \min_{1\le i\le N} |\xi-a_i|^q \mu(d\xi) = 0$. In particular $ \lim_N e_p(a^{(N)},X)=0$.
\end{Pro}
\noindent {\bf Proof.} $(a)$ We proceed by induction. When $N=1$, the existence of $a_1$ is obvious once noticed that $\xi\mapsto \E\, |X-\xi|^p$ is continuous and goes to infinity as $|\xi|\to +\infty$. Assume  there exists $a_1,\ldots a_{_N}$ such that $e_p(a^{(k)},X)= \min_{a\in \R^d} e_p(a^{(k-1)}\cup\{a\},X)$ for every $k\!\in \{1,\ldots,N\}$. 

If ${\rm supp}(\mu)\subset \{a_1,\ldots,a_{_N}\}$ then for every $a\!\in \R^d$,   $e_p(a^{(N)}\cup\{a\},X)= e_p(a^{(N)},X)$. Otherwise, let $\xi^*\!\in {\rm supp}(\mu)\setminus\{a_1,\ldots,a_{_N}\}$. It is clear that $|\xi-\xi^*|< d(\xi,a^{(N)})$ on the ball $B\big(\xi^*, \frac 14d(\xi^*, a^{(N)})\big)$ which satisfies $\mu\big(B\big(\xi^*, \frac 14d(\xi^*, a^{(N)})\big)\big)>0$. Consequently, $e_p(a^{(N)}\cup\{\xi^*\},X) < e_p(a^{(N)},X)$. Now let
\[
K^0_{N+1} = \big\{\xi\!\in \R^d\,|\, e_p(a^{(N)}\cup\{\xi\},X)  \le e_p(a^{(N)}\cup\{\xi^*\},X) \big\}.
\]
This is a closed non-empty set. Now let $(\xi_k)_{k\ge 1}$ be a sequence of elements of $K^0_{N+1}$ such that $|\xi_k|\to +\infty$. It follows from Fatou's Lemma that
\begin{eqnarray*}
\liminf_k e_p(a^{(N)}\cup\{\xi_k\})^p &\ge &\int_{\R^d} \liminf_k \big(d(\xi, a^{(N)})^p\wedge |\xi-\xi_k|^p\big) \mu(d\xi)\\
&=& \int_{\R^d} d(\xi, a^{(N)})^p  \mu(d\xi) \\
&=&e_p(a^{(N)},X)^p> e_p(a^{(N)}\cup\{\xi^*\},X)^p.
\end{eqnarray*}
This yields a contradiction which in turn implies  that $K^0_{N+1}$ is a compact set. On the other hand $\xi\mapsto e_p(a^{(N)}\cup\{\xi\},X)$ is clearly Lipschitz continuous on $\R^d$, hence it attains its minimum on $K^0_{N+1}$ which is clearly its absolute minimum. 

\smallskip
\noindent $(b)$ If $\mu\big(C_{a_N}(a^{(N)})\big)=0$, then one checks that 
$$
e_p(a^{(N-1)}, X)-e_p(a^{(N)},X)= \int_{C_{a_N}(a^{(N)})} \big(d(\xi,a^{(N-1)})^p-|x-a_{_N}|^p\big)\mu(d\xi)=0
$$ 
which contradicts the strict decreasing monotony of $e_p(a^{(N)},X)$.  Let $(C_i)_{1\le i\le N}$ be a Borel Voronoi partition of $\R^d$ induced by $a^{(N)}$, $i.e.$ satisfying $C_i\subset W_{a_i}(a^{(N)})= \big\{\xi\!\in \R^d\,|\, |\xi-a_i| = \min_{1\le j\le N} |\xi-a_j|\big\}$, and such that $C_{_N}=C$. Assume there exists $b\!\in C$ such that $\int_{C} |\xi-a_{N}|^p\mu(d\xi)>\int_{C} |\xi-b|^p\mu(d\xi)$. Then
\begin{eqnarray*}
e_p(a^{(N)},X)^p&=& \sum_{i=1}^{N-1}\int_{C_i} |\xi-a_i|^p\mu(d\xi) +\int_{C}|\xi-a_{_N}|^p\mu(d\xi)\\
&\ge & \sum_{i=1}^{N-1}\int_{C_i} d\big(\xi, a^{(N-1)}\cup\{b\}\big)^p\mu(d\xi) +\int_{C}|\xi-a_{_N}|^p\mu(d\xi)\\
&>&  \sum_{i=1}^{N-1}\int_{C_i} d\big(\xi, a^{(N-1)}\cup\{b\}\big)^p\mu(d\xi) +\int_{C} |\xi-b|^p\mu(d\xi)\quad\mbox{ since}\; \mu(C)>0\\
&\ge&  \sum_{i=1}^{N-1}\int_{C_i} d\big(\xi, a^{(N-1)}\cup\{b\}\big)^p\mu(d\xi) +\int_{C}d\big(\xi, a^{(N-1)}\cup\{b\}\big)^p\mu(d\xi)\\
&=& e_p(a^{(N-1)}\cup\{b\},X)^p
\end{eqnarray*}
which contradicts the minimality of $a_{_N}$.

\smallskip
\noindent $(c)$ Let $p\!\in (0,+\infty)$. It is clear that, for every $\xi\!\in \R^d$, $\min_{1\le i\le N} |\xi-a_i|$ is non-increasing and  converges  toward $\displaystyle \inf_{N\ge 1} |\xi-a_{_N}|$ so that by the monotone convergence theorem, one has
\[
e_p (a^{(N)},X)^p \downarrow \ell_{\infty}:=  \int_{\R^d} \inf_{i\ge 1} |\xi-a_i|^p  \mu(d\xi).
\]
Let $a^{(\infty)}=\{a_{N}, \, N\ge 1\}$. If $\ell_{\infty}\neq 0$, then there exists  $\xi_0\!\in {\rm supp}(\mu)$ such that $\varepsilon_0= d(\xi_0, a^{(\infty)}) >0$. Then, for every $\xi \in B(\xi_0, \frac{\varepsilon_0}{4})$, $d(\xi, a^{(\infty)}) \ge \frac{3}{4}\varepsilon_0$ so that 
\[
\int_{B(\xi_0, \frac{\varepsilon_0}{4})}d\big(\xi, a^{(\infty)}\big)^p \mu(d\xi) \ge \eta_0\quad \mbox{ with }\quad \eta_0= \left(\frac{3\varepsilon_0}{4 }\right)^p \mu\Big(B\big(\xi_0, \frac{\varepsilon_0}{4}\big)\Big).
\]
Now, let $N_0$ be a positive  integer such that, 
\[
\int_{\R^d} d\big(\xi, a^{(N_0)}\big)^p \mu(d\xi)\le \ell_{\infty} +  \frac{\eta_0}{2}\Big(1-\frac{1}{3^p }\Big). 
\]
We consider the $(N_0+1)$-quantizer $a^{(N_0)}\cup\{\xi_0\}$. On the one hand,
\[
\int_{B(\xi_0, \frac{\varepsilon_0}{4})}d\big(\xi, \{a_1,\ldots,a_{N_0},\xi_0\}\big)^p \mu(d\xi) \le \left(\frac{\varepsilon_0}{4}\right)^p \mu\Big(B\big(\xi_0, \frac{\varepsilon_0}{4}\big)\Big)=  \frac{\eta_0}{3^p }
\]
and, on the other hand,  \begin{eqnarray*}
\int_{B(\xi_0, \frac{\varepsilon_0}{4})^c}d\big(\xi,  a^{(N_0)}\cup\{
\xi_0\}\big)^p \mu(d\xi) &\le&   \int_{B(\xi_0, \frac{\varepsilon_0}{4})^c}d\big(\xi, a^{(N_0)}\big)^p \mu(d\xi)   \\
& \le  & \int_{\R^d}d\big(\xi, a^{(N_0)}\big)^p \mu(d\xi) -\int_{B(\xi_0, \frac{\varepsilon_0}{4})}d\big(\xi,a^{(N_0)}\big)^p \mu(d\xi) \\
&\le & \ell_{\infty} +  \frac{\eta_0}{2}\Big(1-\frac{1}{3^p }\Big) -  \eta_0
\end{eqnarray*} 
so that
\[
\int_{\R^d}d\big(\xi,  a^{(N_0)}\cup\{
\xi_0\}\big)^p \mu(d\xi) \le \ell_{\infty} +\frac{\eta_0}{2}\Big(1-\frac{1}{3^p }\Big)-  \eta_0 + \frac{\eta_0}{3^p }<\ell_{\infty}
\]
which yields a contradiction. Hence $\ell_{\infty}=0$ which completes the proof for $q=p$. 

\smallskip 
Finally, we derive from what precedes    that $\lim_N \min_{1\le i\le N} |X-a_i| = \inf_{N\ge 1} |X-a_{_N}| =0$ $\P$-$a.s.$. As $X\!\in L^q(\P)$, $ \min_{1\le i\le N} |X-a_i|\le|X-a_1|\!\in L^q$, the conclusion follows from the Lebesgue dominated convergence theorem. $\quad\cqfd$

\bigskip
\noindent {\bf Remark on uniqueness.} Uniqueness of $L^p$-optimal  greedy quantization sequence turns out to be quite different problem from its counterpart for regular $L^p$-optimal quantization. Thus,  for $1$-dimensional  $\log$-concave distributions, it is well-known that uniqueness of $L^p$-optimal quantizers holds true (up to a reordering of the components in an increasing order,  see~\cite{KIE}) holds true. For $L^p$-optimal greedy quantization,  this uniqueness may fail. Basically, greedy quantization is more influenced by the symmetry properties of the distributions: thus for the ${\cal N}(0;1)$-distribution (whose density is $\log$-concave), its is clear that $a_1=0$ (unique $L^p$-median) but then we have that, if $a_2$ is  the (unique, see Proposition~\ref{pro:lloydI-I} in Appendix~\ref{app:A}) solution to the the problem
\[
\min_{a\ge 0}\E \big(|X|^p\wedge|X-a|^p\big)\quad \mbox{ where $X$ has distribution } \mu = {\cal N}(0;1),
\]
 then both $a_2$ and $-a_2$  are solutions to the greedy problem~\eqref{eq:greedy} at level $N=2$ by symmetry of  (the distribution of) $X$. In fact, one derives in turn that   $(0,a_2,-a_2)$ and $(0, -a_2, a_2)$ are both the first three terms of $(L^p,{\cal N}(0;1))$-optimal greedy quantization sequences.

\subsection{About $L^p$-optimal greedy quantization in an Euclidean framework}

In this section we assume that $|\,.\,|$ denotes an Euclidean norm on $\R^d$. Let $\overline {\cal H}_{\mu}$ be the {\em closed} convex hull of the support of the distribution $\mu$.

\begin{Pro}\label{pro:geogreedy2} 
Let $(\,.|.\,)$ denote the  inner product induced by the Euclidean norm.  If ${\rm supp}(\mu)$ contains at least $N$ elements then, the  first $N$ elements of any optimal greedy quantization sequence takes values in $\overline{\cal H}_{\mu}$. If ${\rm supp}(\mu)$ is infinite any optimal greedy quantization sequence takes values in $\overline{\cal H}_{\mu}$.  
\end{Pro}

\noindent {\bf Proof.} We proceed by induction. Let $\displaystyle a_1\!\in {\rm argmin}_{a\in \R^d}\E|X-a|^p$ and let $\Pi_1(a_1)$ be  its projection on $\overline {\cal H}_{\mu}$. If $a_1\neq \pi_1(a_1)$, the pseudo-Pythagoras Theorem implies
\[
\forall\, \xi\!\in \overline {\cal H}_{\mu},\quad |\xi-a_1|^2\ge |\xi-\pi_1(a_1)|^2 + |a_1-\pi_1(a_1)|^2
\]
so that $\E\,|X-a|^2\ge \E|X-\pi_1(a_1)|^2 +  |a_1-\pi_1(a_1)|^2$ which yields a contradiction to the definition of $a_1$. Hence $a_1\!\in  \overline {\cal H}_{\mu}$.

\smallskip Let $a_{_N}\!\in {\rm argmin}_{a\in \R^d}e_p(a^{(N-1)}\cup\{a\},X)$. It follows from Proposition~\ref{pro:greedex}$(b)$ that $a_{_N}\!\in {\rm argmin}_{a\in \R^d}\E\big(|X-a|^p\,|\, X\!\in W_{a_N}\big)$ where $W_{a_N}= \{\xi\!\in \R^d\,|\, |\xi-a_{_N}|\le d(\xi,a^{(N-1)})\}$ is a  closed (polyhedral) convex set since the norm is Euclidean and has a positive $\mu$-measure. As a consequence $a_{_N}\!\in \overline{\cal H}_{\mu(.\,|W_{a_N})}$ where $\mu(.\,|W_{a_N})$ is the conditional distribution of $\mu$ given $W_{a_N}$. One concludes by noting that $ \overline{\cal H}_{\mu(.\,|W_{a_N})}= \overline{\cal H}_{\mu} \cap W_{a_N}\subset  \overline{\cal H}_{\mu}$.~$\cqfd$

%
%

\bigskip
\noindent {\bf Remark.} Let $p=2$. As soon as ${\rm card}({\rm supp}(\mu))\ge N$, we know from Proposition~\ref{pro:greedex}$(b)$  that  $\mu(W_{a_N})>0$  and
\[
{\rm argmin}_{a\in \R^d} \int_{W_{a_N}}|\xi-a|^2\mu(d\xi) =\left\{\frac{ \int_{W_{a_N}}\!\xi \,\mu(d\xi) }{\mu(W_{a_N})}\right\}
\]
$i.e.$ 
\begin{equation}\label{eq:Lloydloc}
a_{_N}= \frac{ \int_{W_{a_N}}\!\xi \,\mu(d\xi) }{\mu(W_{a_N})}=\E\big(X\,|\,X\in W_{a_N} \big).
\end{equation}
 This can be seen as a fixed point formula and is the starting point of stochastic optimization procedure to compute by simulation (of i.i.d. samples of $X$) of optimal greedy sequences  using a variant of the celebrated Lloyd method introduced in~\cite{LLOYD} and widely used in Statistics and Data Analysis(see~\cite{McQ})  as $k$-means algorithm  (see Section~\ref{sec:greedyalgo}).
%
%

\section{Greedy quantization is rate optimal}\label{sec:Rate}
\subsection{A general rate optimality result}
Following~\cite{GrLuPa2}, we define for every $b\!\in (0,\frac12)$ the {\em  $b$-maximal function} associated to an $L^p$-optimal  greedy quantization sequence $(a_N)_{N\ge 1}$  by
\[
\forall\, \xi\!\in \R^d,\quad \Psi_{b}(\xi) = \sup_{N\ge 1} \frac{\lambda_d\big(B(\xi,bd(\xi,a^{(N)}))\big)}{\mu\big(B(\xi,bd(\xi,a^{(N)}))\big)}\!\in [0,+\infty].
\]
It is clear that $\Psi_b(\xi) >0$ for every $\xi\neq a_1$ ($L^p$-median).  

\smallskip
 Note that this notion of $b$-maximal function (originally introduced in~\cite{GrLuPa2}) can be naturally defined with respect to a sequence of grids $(\Gamma_{N})_{N\ge 1}$ where $\Gamma_{N}$ has size $N$. 
 
 The theorem below yields a criterion based on the integrability of the maximal function $\Psi_b$
which implies that an $(L^p,\mu)$-optimal  greedy quantization sequence is $(L^p,\mu)$-rate optimal (in the sense of Zador's Theorem). More practical criteria are given further on in Section~\ref{sec:PractiCriterion}. 

\begin{Thm}\label{thm:main} Let $p\!\in (0,+\infty)$ and let $\mu=\P_{_X}$ be such that    $\displaystyle \int_{\R^d}|\xi|^p\mu(d\xi)<+\infty$.   Let $(a_N)_{N\ge 1}$ be an $L^p$-optimal  greedy quantization sequence. Assume that there exists $b\!\in (0, \frac 12)$ such that $\Psi_b \!\in L^{\frac{p}{p+d}}(\mu)$. Then
\begin{equation}\label{eq:RateOpt}
\limsup_N  N^{\frac 1d} e_p(a^{(N)}, X) <+\infty.
\end{equation}
\end{Thm}

\noindent {\bf Proof.} First, note that if $\mu$ is a Dirac mass $\delta_a$ for some  $a\!\in \R^d$, then $a_1=a$ and $e_p(a^{(N)},X)=0$ for every integer $N\ge1$.  Otherwise, we rely on the following  {micro-macro inequality} established in~\cite{GrLuPa2} (see Equation~(3.4) in the proof of Theorem~2, with the standard convention $\frac 10=+\infty$).
\[
\forall\, \xi\!\in \R^d,\qquad d(\xi,a^{(N)})^p \le \frac{C_{p,b}}{\mu\big(B(\xi,bd(\xi,a^{(N)}))\big) }\Big(e_p(a^{(N)},X)^p - e_p(a^{(N)}\cup\{\xi\},X)^p\Big)
\]
where $b\!\in (0,\frac12)$ and $C_{p,b}$ is a positive real constant depending on $p$ and $b$. Then, it follows that
\begin{equation}\label{eq:micro-macro}
e_p(a^{(N)}\cup\{\xi\},X)^p\le e_p(a^{(N)},X)^p -\frac{1}{C_{p,b}} \frac{\mu\big(B(\xi,bd(\xi,a^{(N)}))\big)}{\lambda_d\big(B(\xi,bd(\xi,a^{(N)}))\big)}b^d d(\xi, a^{(N)} )^{p+d}V_d
\end{equation}
where $V_d$ denotes the hyper-volume of the  unit ball with respect to the current norm on $\R^d$ $i.e.$ $V_d = \lambda_d\big(B_{|\,.\,|}(0;1)\big)$. This implies that
\begin{equation}\label{eq:micro-macro2}
e_p(a^{(N)}\cup\{\xi\},X)^p\le e_p(a^{(N)},X)^p -\frac{1}{\widetilde C_{p,b,d}} \frac{1}{\Psi_{b}(\xi)} d(\xi, a^{(N)})^{p+d}
\end{equation}
where $\widetilde C_{p,b,d}= C_{p,b}/( b^d V_d)\!\in(0,+\infty)$. Note that $\mu(\{a_1\})<1$  since $\mu$ is not a Dirac mass, so that 
\[
\int_{\R^d} \Psi^{\frac{p}{p+d}}_b(\xi)\,\mu(d\xi) >0.
\]
Consequently, as $\Psi_b \!\in L^{\frac{p}{p+d}}(\mu)$, we can define the probability distribution $\displaystyle  \nu= \kappa_{b,p,d}\,\Psi_b^{\frac{p}{p+d}}.\mu$ (where $\kappa_{b,p,d}= \Big(\int\Psi_b^{\frac{p}{p+d}}d\mu\Big)^{-1}\!\in (0,+\infty)$ is a normalizing real constant). Then, integrating the above inequality with respect to $\nu$ yields
\[
\int_{\R^d}e_p\big(a^{(N)}\cup\{\xi\}, X\big)^p \nu(d\xi)  \le e_p(a^{(N)},X)^p - \widetilde C_{p,b,d}\int_{\R^d} d(\xi,a^{(N)})^{p+d} \frac{\nu(d\xi)}{\Psi_b(\xi)}.
\]
Jensen's Inequality applied  to the convex function $u\mapsto u^{1+\frac dp}$ yields  
\begin{eqnarray*}
\int_{\R^d} d(\xi,a^{(N)})^{p+d} \frac{\nu(d\xi)}{\Psi_b(\xi)}& \ge & \left(\int_{\R^d} d(\xi,a^{(N)})^{p} \frac{\nu(d\xi)}{\Psi_b(\xi)^{\frac{p}{p+d}}}\right)^{1+\frac dp}\\
&=& \kappa_{b,p,d}^{1+\frac dp} \left(\int_{\R^d} d(\xi,a^{(N)})^{p} \mu(d\xi)\right)^{1+\frac dp}\\
&=&  \kappa_{b,p,d}^{1+\frac dp} \,e_p\big(a^{(N)},X\big)^{p+d}.
\end{eqnarray*}
On the other hand, it is clear that
\[
e_p\big(a^{(N+1)},X\big)^p\le \int_{\R^d}\nu(d\xi) e_p\big(a^{(N)}\cup\{\xi\}, X)\big)^p
\]
so that, finally, if we set $A_{_N} = e_p(a^{(N)},X)^p$, $N\ge 1$, this sequence satisfies for every integer $N\ge 1$,  the recursive inequality
\[
A_{_{N+1}}\le A_{_N}- \widetilde \kappa\, A_{_N}^{1+\frac d p}
\]
where $\widetilde \kappa = \kappa_{b,p,d}^{1+\frac dp}\widetilde C_{p,b,d}$. The sequence $(A_{_N})_{N\ge 1}$ being non-negative,  one classically derives the announced conclusion (for a proof, see Lemma~\ref{lem:A-N} in   the Appendix~\ref{app:B}, applied with $\rho = \frac dp$ and $C=\widetilde \kappa$).~$\cqfd$

\bigskip
\noindent {\bf Remark.}  $\bullet$ One straightforward derives from Zador's Theorem (Theorem~\ref{thm:Zador}$(a)$) that,  under the assumption of the above theorem 
and if $\mu$ has a non-zero absolutely continuous component ($i.e.$ $\varphi=\frac{d\mu}{d\lambda_d}  \not\equiv 0$), one has
\[
e_p(a^{(N),p},X)\asymp N^{-\frac 1d}
\]
since $e_p(a^{(N),p},X)\ge e_{p,N}(X)$ and  $\displaystyle \liminf_N N^{\frac 1d} e_{p,N}(X) \ge \widetilde J_{p,d} \|\varphi\|_{L^{\frac{p}{p+d}}(\lambda_d)}^{\frac{1}{p}} >0$.
The same conclusion will hold true for the distortion mismatch problem investigated in Proposition~\ref{pro:greedymismatch} in the next section.

\smallskip
\noindent $\bullet$ A careful reading of the proof shows that, if we define the sequence of functions $\Psi_{b,N}$ by 
\[
\forall\, \xi\!\in \R^d,\quad \Psi_{b,N}(\xi) = \frac{\lambda_d\big(B(\xi,bd(\xi,a^{(N)}))\big)}{\mu\big(B(\xi,bd(\xi,a^{(N)}))\big)}\!\in [0,+\infty],
\]
then the theorem  holds true under the weaker assumption that  there exists an integer $N_0\ge 1$ such that $\displaystyle \sup_{N\ge N_0}\int_{\R^d}\Psi_{b,N}(\xi)\mu(d\xi)<+\infty$.Unfortunately, this fact seems to be of little practical interest. 

\smallskip
\noindent  $\bullet$ When $\mu$ is singular with respect to the Lebesgue measure (no absolutely continuous part), it is likely that, like for standard optimal vector quantization in Zador's Theorem, this rate is not optimal. The natural conjecture should be that  greedy quantization sequence(s) go to $0$ at the same rate as that obtained for sequences of optimal quantizers which is not $N^{-\frac 1d}$ when the distribution $\mu$ is singular (see $e.g.$~\cite{GRLU}).

\smallskip 
\noindent $\bullet$   Since we know that $d(\xi, a^{(N)})\downarrow 0$ as $N\to +\infty$, $\mu(d\xi)$-$a.s.$, it is  clear that  if 
$\mu =\varphi.\lambda_d\displaystyle$ (or even $\mu =\varphi.\lambda_d\displaystyle  \stackrel{\perp}{+} \tilde \mu$, to be checked), then 
by the Lebesgue  differentiation theorem
\[
\frac{1}{\varphi(\xi)}= \liminf_N \frac{\lambda_d\big(B(\xi,bd(\xi,a^{(N)}))\big)}{\mu\big(B(\xi,bd(\xi,a^{(N)}))\big)}\le \Psi_b(\xi) \quad\mu(d\xi)\mbox{-}a.s.
\]
so that by Fatou's Lemma, the condition $\Psi_b \!\in L^{\frac {p}{p+d}}(\mu)$ implies 
\[
\int_{\R^d} \varphi^{\frac{d}{p+d}}(\xi)d\lambda_d(d\xi)<+\infty.
\]
So, we retrieve here the statement of Remark~6.3$(c)$, p.79,  in~\cite{GRLU} which points out  that  {\em if} optimal  $L^p$-mean quantization goes to zero at rate $N^{-\frac 1d}$ then the above integral is finite (see also Section~1 in~\cite{GrLuPa2}). Of course, as emphasized in Remark~6.3$(a)$ from~\cite{GRLU}, p.79,  the classical condition under which Zador's Theorem holds, namely  $\E|X|^{p+\delta}=\int_{\R^d}|\xi|^{p+\delta}\mu(d\xi)<+\infty$ for a $\delta>0$, implies the finiteness of this integral owing to an appropriate application of H\"older's inequality. The above result suggests a hopefully nonempty question:   since $L^p$-rate optimality for greedy sequence (and consequently for true $L^p$-optimal quantizers) holds as soon as $X\!\in L^p(\P)$ and $\psi_b(X)\!\in L^{\frac{p}{p+d}}(\P)$ for a $b\!\in (0,\frac 12)$, are such conditions achievable when $\E |X|^{p+\delta}=+\infty$ for every $\delta>0$.

 \subsection{Distortion mismatch for optimal greedy quantization sequences}
 
 In this section we address the problem of distortion mismatch originally investigated in~\cite{GrLuPa2} for sequences of optimal $N$-quantizers.  
 
 If $q\!\in(0,p]$ and $X\!\in   L^p(\P)$ any optimal greedy sequence $(a_{_N})_{N\ge 1}$ remains $L^q$-rate optimal for the $L^q$-norm owing to the monotonicity of the $L^q$-norm as function of $q$. But the challenging question for distortion mismatch starts with the case $q>p$. It is solved in the proposition below, still relying on an integrability assumption on the $b$-maximal function(s) $\Psi_b$. For more practical criteria we again  refer to Section~\ref{sec:PractiCriterion}. 
 
 \begin{Pro} \label{pro:greedymismatch}Let $q\!\in (p,+\infty)$ and let $X\!\in L^p(\P)$ with distribution $\mu=\P_{_X}$. Assume that 
 the maximal function $\Psi_b \!\in L^{\frac{q}{p+d}}(\mu)$ for some $b\!\in (0, \frac 12)$.  Let $(a_N)_{N\ge 1}$ be an $L^p$-optimal greedy sequence.
 
 Then  $X\!\in L^q(\P)$ and 
 \[
\limsup_N N^{\frac 1d} e_{q}(a^{(N)},X) <+\infty.
 \]
%
\end{Pro}

\noindent {\bf Remarks.}  When ${\rm supp}(\mu)$ is not compact it is hopeless  to have results for $q>p+d$  since it has been shown in~\cite{GrLuPa2} (Theorem~10 and Equation (2.7)) that  the $L^q$-rate optimality of a sequence $(a_N)_{N\ge 1}$ would imply when $\mu=\varphi.\lambda_d$ that 
$$
\int_{\varphi>0}\varphi^{-\frac{q}{p+d}}(\xi)\mu(d\xi)= \int_{\varphi>0}\varphi^{1-\frac{q}{p+d}}(\xi)\lambda_d(d\xi)<+\infty.
$$
However when $\mu$ has a compact support, we will see in Proposition~\ref{pro:compact2}$(c)$ that $L^q$-rate optimality can be preserved  under appropriate  integrability assumptions.


\bigskip

\noindent {\bf Proof.}  First, note that if $\mu$ is a Dirac mass $\delta_a$ for some  $a\!\in \R^d$, then $a_1=a$ and $e_q(a^{(N)},X)=0$ for every integer $N\ge1$. Otherwise,  it follows from Equation~\eqref{eq:micro-macro} rewritten in a reverse  way that 
\[
\forall\, \xi\!\in \R^d, \quad d(\xi,a^{(N)})^q\le C_{b,d,p,q}\left(e_p(a^{(N)},X)^p-e_p(a^{(N)}\cup\{\xi\},X)^p\right)^{\frac{q}{p+d}}\Psi_b(\xi)^{\frac{q}{p+d}}(\xi).
\]
Now,  we note that 
\[
\forall\, \xi\!\in \R^d, \quad e_p(a^{(N)}\cup\{\xi\},X)^p \ge e_p(a^{(N+1)},X)^p
\]
by definition of the sequence $(a_N)_{N\ge 1}$ so that 
\[
\forall\, \xi\!\in \R^d, \quad d(\xi,a^{(N)})^q\le C_{b,d,p,q}\left(e_p(a^{(N)},X)^p-e_p(a^{(N+1)},X)^p\right)^{\frac{q}{p+d}}\Psi_b(\xi)^{\frac{q}{p+d}}(\xi).
\]
Integrating with respect to $\mu$ yields
\[
e_q(a^{(N)},X)^q\le C_{b,d,p,q} \left(e_p(a^{(N)},X)^p-e_p(a^{(N+1)},X)^p\right)^{\frac{q}{p+d}}
\int_{\R^d} \Psi_b(\xi)^{\frac{q}{p+d}}(\xi)\mu(d\xi).
\]
 We know that $\displaystyle \int_{\R^d} \Psi_b(\xi)^{\frac{q}{p+d}}(\xi)\mu(d\xi)\!\in (0,+\infty)$ owing to the assumption made on $\mu$ and $\psi_b$.
Hence 
 \[
e_q(a^{(N)},X)^q\le \widetilde C_{b,d,p,q}   \left(e_p(a^{(N)},X)^p-e_p(a^{(N+1)},X)^p\right)^{\frac{q}{p+d}}\]
where $ \widetilde C_{b,d,p,q}= C_{b,d,p,q} \int_{\R^d} \Psi_b(\xi)^{\frac{q}{p+d}}(\xi)\mu(d\xi)$.
Equivalently
\begin{equation}\label{eq:keyb}
e_q(a^{(N)},X)^{p+d}\le \widetilde C_{b,d,p,q}^{\frac{p+d}{q}}\left(    e_p(a^{(N)},X)^p - e_p(a^{(N+1)},X)^p\right). 
\end{equation}

Summing over $k$ between $N$ and $2N-1$ yields
\begin{eqnarray*}
  \sum_{k=N}^{2N-1}e_q(a^{(k)},X)^{p+d}&\le& \widetilde C_{b,d,p,q}^{\frac{p+d}{q}} \Big( e_p(a^{(N)},X)^p -e_p(a^{(2N)},X)^p\Big)\\
 &\le & \widetilde C_{b,d,p,q}^{\frac{p+d}{q}} 
 e_p(a^{(N)},X)^p.
\end{eqnarray*}

It is clear that $\Psi_b \!\in L^{\frac{p}{p+d}}(\mu)$ since $p<q$ and $\Psi_b \!\in L^{\frac{q}{p+d}}(\mu)$. Consequently, 
 it follows from 
 Theorem~\ref{thm:main} 
 that there exists a positive real constant $\widetilde C'_{b,d,p,q}\!\in (0, +\infty)$ such that, for every $N\ge 1$, 
 \[
  \sum_{k=N}^{2N-1}e_q(a^{(k)},X)^{p+d}\le  \widetilde C'_{b,d,p,q}N^{-\frac pd}.
 \]
 On the other hand the sequence $\big(e_q(a^{(N)}, X)\big)_{N\ge 1}$ is clearly non-decreasing since $(d(\xi,a^{(N)})^q)_{N\ge 1}$ is itself non-decreasing  for every $\xi\!\in \R^d$. Finally, this implies that, for every $N\ge1$,
 \[
Ne_q(a^{(2N-1)},X)^{p+d}\le  \sum_{k=N}^{2N-1}e_q(a^{(k)},X)^{p+d}\le  \widetilde C'_{b,d,p,q}N^{-\frac pd}.
\]
Hence, for every integer $N\ge 1$, 
\[
Ne_q(a^{(N)},X)^{p+d}\le 2 \,  \widetilde C'_{b,d,p,q}\lceil N/2\rceil^{-\frac pd}.
\]
(where $\lceil m \rceil$ denotes the upper integer part of $m\!\in\N$). 
Consequently,  for every $N\ge1$,
\[
e_q(a^{(N)},X)^{p+d}\le2^{1+\frac pd} \widetilde C'_{b,d,p,q}N^{-(1+\frac pd)}.
\]
One completes the proof by taking the $(n+p)^{th}$ root of the inequality.~$\cqfd$

 \section{Practical criteria for the integrability of the  maximal function} \label{sec:PractiCriterion}

 These criteria are mainly borrowed from~\cite{GrLuPa2} where they have   been established for the first time in order to solve the mismatch problem for optimal quantization 

\paragraph{Compact case and $q<p+d$.} The compact case relies on the following lemma which allows for non convex support for the distribution $\mu$.
\begin{Lem}[see Lemma~1 in~\cite{GrLuPa2}] If $X\!\in L^p(\P)$ has  a distribution $\mu$ 
and  $(\Gamma_N)_{N\ge 1}$ is a sequence of $N$-quantizers such that $\int_{\R^d} d(\xi,\Gamma_{_N})^p\mu(d\xi)\to 0$, then the maximal functions $\Psi_b$ associated to $(\Gamma_{_N})_{N\ge 1}$  lie in $L^r_{loc}(\mu)$ for every  $r\!\in (0,1)$ $i.e.$
\[
\forall\, r\!\in (0,1),\; \forall\, b\!\in (0,\frac 12),\;\forall R\!\in (0,+\infty), \quad \int_{\{|\xi|\le R\}}\psi_b(\xi)^r \mu(d\xi)<+\infty.
\]
\end{Lem}

By combining this result (applied with $r= \frac{q}{p+d}$) with Proposition~\ref{pro:greedex}$(b)$, we derive the following result  which extends the one established in~\cite{Brancoetal} for absolutely continuous distributions with convex support on $\R^d$.  Note that the proof of the above lemma is not elementary, especially when ${\rm supp}(\mu)$ is not convex,  and relies on the Besicovitch covering theorem.

\begin{Pro}[Compact support]  If $X$ has a distribution $\mu$ with compact support, then any $L^p$-optimal greedy quantization sequence $(a_{_N})_{N\ge 1 }$ is $L^q$-rate optimal for every $q\!\in (0,p+d)$ $i.e.$ satisfies 
\[
\limsup_N N^{\frac1d} e_q(X,a^{(N)})<+\infty.
\]
\end{Pro}

\paragraph{Compact case and $q\ge p+d$.} 
Results can be derived for $q>p+d$  when $\mu$ is absolutely continuous and has a compact support. They   rely on the following Lemma (see Lemma~2  in~\cite{GrLuPa2}).
\begin{Lem} Assume $\mu= \varphi.\lambda_d$, $\E |X|^p<+\infty$, ${\rm supp}(\mu)$ is the finite union of closed convex sets and $\lambda_{d|{\rm supp}(\mu)}$ is absolutely continuous with respect to $\mu$.

 Let $(\Gamma_{_N})_{N\ge 1}$ be a sequence of quantization grids satisfying $e_p(\Gamma_{_N}, X)\to 0$ as $N\to +\infty$.  Then, for every $q\!\in (1,+\infty]$, the associated maximal functions $\Psi_b$ lie in $L^q_{loc}(\mu)$ iff $\frac{1}{\varphi}\!\in L_{loc}^q(\mu)$. 
\end{Lem}

As a consequence of this lemma, we derive the following proposition which deals with the cases $q>p+d$ (in $(a)$) and $q=p+d$ (in~$(b)$).

\begin{Pro}\label{pro:compact2}  $(a)$  Let  $\mu= \varphi.\lambda_d$ be like in the preceding lemma and let $(a_N)_{N\ge 1}$ be  an $L^p$-optimal greedy quantization  sequence for $\mu$. Let $q>d+p$. If
\[
\int_{\R^d} \varphi^{-\frac{q}{d+p}}(\xi)\mu(d\xi) = \int_{\{\varphi>0\}} \varphi^{1-\frac{q}{d+p}}(\xi)\lambda_d(d\xi) <+\infty
\]
then $(a_N)_{N\ge 1}$ is $L^{q'}$-rate optimal for every $q'\!\in (0,q]$ $i.e.$
\[
\limsup_N N^{\frac 1d} e_{q'}(X, a^{(N)})<+\infty.
\]
In particular,  if  $\varphi\ge \varepsilon>0$ on ${\rm supp}(\mu)$, then the above integral criterion   is fulfilled.

\smallskip
\noindent $(b)$  Let $q=p+d$. If there exists $\delta>0$ such that 
\[
 \int_{\R^d} \varphi^{-(1+\delta)}(\xi)\mu(\xi)=\int_{\{\varphi>0\}} \varphi^{-\delta}(\xi)\lambda_d(d\xi) <+\infty
\]
 then $\hskip 4 cm \displaystyle \limsup_N N^{\frac 1d} e_{p+d}(X, a^{(N)})<+\infty$.
\end{Pro} 

\paragraph{Non-compact radial case.}
\begin{Lem}[see Corollary 3 in~\cite{GrLuPa2}]  If $X\!\in L^{p+\delta}(\P)$ for some $\delta>0$  with an essentially radial distribution $\mu(d\xi)=\varphi(\xi)\lambda_d(d\xi)$ in the sense  that
\begin{equation}\label{eq:radial}
\hskip -0,25cm \varphi=h(|\,.\,|_0)\, \mbox{ on }\; B_{|\,.\,|_0}(0,R)^c \mbox{ with } h:(R,+\infty)\to \R_+,\mbox{ non-increasing and $|\,.\,|_0$ any norm on $\R^d$.}
\end{equation}
Let $(\Gamma_N)_{N\ge 1}$ be a sequence of $N$-quantizers such that 
$e_q(\Gamma_{_N}, X) \to0$. If there exists a real constant $c>1$ such that 
\begin{equation}\label{eq:fc}
\int_{\R^d}\varphi(c\,\xi)^{-\frac{q}{p+d}}\mu(d\xi)= \int_{\R^d}\varphi(c\,\xi)^{-\frac{q}{p+d}}\varphi(\xi)d\xi <+\infty
\end{equation}
then  $\Psi_b\!\in L^{\frac{q}{p+d}}(\mu)$.
\end{Lem}

In fact, as stated in~\cite{GrLuPa2}, Corollary~3 is written to be used only with $L^p$-optimal quantizers so the above formulation includes minor modifications. Combining this lemma with Proposition~\ref{pro:greedex}$(b)$ and Theorem~\ref{thm:main} yields the following proposition.

\begin{Pro}[Non-compact support with radial density]  If $X\!\in L^{p+\delta}(\P)$ for some $\delta>0$ with an essentially  radial distribution in the sense of~\eqref{eq:radial} and if, furthermore, $\varphi$ satisfies~\eqref{eq:fc}, then any  $L^p$-optimal greedy sequence $(a_{_N})_{N\ge 1 }$ is $L^q$-rate optimal $i.e.$ satisfies 
\[
\limsup_N N^{\frac1d} e_q(X,a^{(N)})<+\infty.
\]
\end{Pro}
This case includes $e.g.$ all the  centered {\em hyper-exponential distributions}   of the form $\mu= \varphi.\lambda_d$ with 
\[
\varphi(\xi)= \kappa_{a,b,c}|\xi|_0^c \, e^{-a|\xi|_0^b},\; \xi\!\in \R^d,\; a,b>0,\; c>-d
\]
and $|\,.\,|_0$ is any norm on $\R^d$ and subsequently  all hyper-exponential distributions since $L^p$-mean-quantization errors are invariant by translation of the random vector $X$. In particular, this  includes  all   normal and Laplace distributions.

\bigskip
\noindent {\bf Remark.} In one dimension,~\eqref{eq:radial} can be replaced {\em mutatis mutandis} by  a one-sided variant: if there exist $R_0, R'_0\!\in \R$,  $R'_0\ge R_0$ such that
\begin{equation}\label{eq:radial1}
{\rm supp}(\mu)\subset [R_0,+\infty)\;\mbox{ and }\;   f_{|[R'_0,+\infty)}\mbox{ is non-increasing}.
\end{equation}

This criterion is satisfied by the gamma distributions on $\R_+$ (including the exponential distributions).


\paragraph{Non-compact and possibly non-radial case.}
\begin{Cor}\label{Cor4} Assume $\mu =  \varphi. \lambda_d$ and $\displaystyle\E\,|X|^{p+\delta} <+ \infty$ for some $\delta > 0$. Furthermore, assume that ${\rm supp} (\mu)$ is {\em peakless}  in the following sense
\begin{equation}\label{peak}
\kappa_{\varphi}:=\inf_{\xi\in {\rm supp} (\mu),\, 0<\rho\le 1} \frac{\lambda_d \big({\rm supp} (\mu)\cap B(\xi,\rho)\big)}{\lambda_d \big(B(\xi,\rho)\big)}>0
\end{equation}
and that $\varphi$ satisfies the {\em  local growth control assumption}: there exist    real numbers
$\varepsilon \ge 0$, $\n\!\in(0,\frac 12)$, $M,\,C>0$ such that
\begin{equation}\label{Cond2}
 \forall\, \xi,\,\xi' \!\in{\rm supp }(\mu),\; |\xi|\ge M,\;|\xi'-\xi|\le
2\n\, |\xi|\;\Longrightarrow \;\varphi(\xi') \ge C \varphi(\xi)^{1+\varepsilon}.
\end{equation}
Then, for every   $q\!\in(0,\frac{p+d}{1+\varepsilon})$ such that
$$
 \int_{\R^d}\varphi(\xi)^{-\frac{q(1+\varepsilon)}{p+d}}\mu(d\xi)= \int_{\{\varphi (\xi)>0\}}
 \varphi(\xi)^{1-\frac{q(1+\varepsilon)}{p+d}}\lambda_d(d\xi) <+\infty
$$
 (if any),  any greedy $L^p$-optimal sequence $(a_{_N})_{N\ge 1 }$ is $L^q$-rate optimal $i.e.$ satisfies 
\[
\limsup_N N^{\frac1d} e_q(X,a^{(N)})<+\infty.
\]

In particular, if~(\ref{Cond2}) holds either for  $\varepsilon=0$ or for every $\varepsilon \!\in(0,\underline
\varepsilon]$ $(\underline \varepsilon>0$), and if
\begin{equation}\label{CondInteg}
\forall\, q\!\in(0,p+d),\qquad  \int_{\R^d}
\varphi(\xi)^{-\frac{q}{p+d}}\mu(d\xi)=\int_{\{\varphi (\xi)>0\}}
 \varphi(\xi)^{1-\frac{q}{d+p}}\lambda_d(d\xi)  <+\infty
\end{equation}
then the above conclusion  holds   for every $ q\!\in(p,p+d)$.
\end{Cor}

Note that (if $\lambda_d({\rm supp}(\mu))=+\infty$)
Assumption~(\ref{peak})   is $e.g.$ satisfied by any finite intersection of half-spaces, the typical
example being
$\R_+^d$. Furthermore, a careful reading of the proof below shows that this assumption
can be slightly relaxed into: there exists a real $c>0$ such that
\[
\kappa'_{\varphi}:=\inf_{\xi\in {\rm supp} (\mu)} \left\{\frac{\lambda_d ({\rm supp}
(\mu)\cap B(\xi,\rho))}{\lambda_d (B(\xi,\rho))},\; 0<\rho \le c\,|x\
\right\}>0.
\]

\section{Further answers  and questions about greedy quantization} \label{sec:Few}

In this section, we temporarily denote by $ \big(a_{N,p}\big)_{N\ge 1}$ the $L^p$-optimal greedy quantization sequence for the uniform distribution $U([0,1])$ and by $\big(\a^{(N),p}\big)_{N\ge 1}$ the resulting  sequence of greedy quantizers.

\paragraph{$\rhd$ Rate optimality of greedy sequences}
It is a straightforward consequence of Zador's Theorem  that if the distribution $\mu=\P_{_X}$ of $X\!\in L^{p+\delta}$, $\delta>0$, has a non-zero absolutely continuous component ($i.e.$ $\varphi=\frac{d\mu}{d\lambda_d}  \not\equiv 0$) and satisfies the assumptions of Theorem~\ref{thm:main}, then   
\[
e_p(a^{(N),p},X)\asymp N^{-\frac 1d}
\]
since $e_p(a^{(N),p},X)\ge e_{p,N}(X)$ and  $\liminf_N N^{\frac 1d} e_{p,N}(X) \ge \widetilde J_{p,d} \|\varphi\|_{L^{\frac{d}{p+d}}(\lambda_d)}^{\frac{1}{p}} >0$. (By the way it proves that under the assumption of Theorem~\ref{thm:main}, $ \|\varphi\|_{L^{\frac{d}{p+d}}(\lambda_d)}<+\infty$.) 

By a similar argument, the same holds true for the distortion mismatch problem under the assumptions of Proposition~\ref{pro:greedymismatch}.

\paragraph{$\rhd$ Can greedy quantization sequence produce asymptotically optimal quantizers?} If $\mu$ has an absolutely continuous component with density $\varphi$,
then  any sequence $(\Gamma_n)_{n\ge 1}$ of {\em asymptotically $(L^p,\mu)$-optimal quantization grids}  at level $N_n={\rm card}(\Gamma_{n})\to +\infty$, satisfies the {\em empirical measure theorem}  (see~\cite{GRLU}, Theorem~7.5, p.96 and~\cite{DeGrLuPa} for a slight refinement), namely 
\[
\frac{1}{N_n} \sum_{a\in \Gamma_n} \delta_a \stackrel{(w)}{\longrightarrow} \mu^{(p)} = \frac{\varphi^{\frac{d}{d+p}}}{\int_{\R^d}\varphi^{\frac{d}{d+p}}d\lambda_d }. \lambda_d\quad\mbox{ as }\quad n\to+\infty
\]
where $\stackrel{(w)}{\longrightarrow}$ denotes the weak convergence of probability measures. Note that when $\mu=U([0,1])$, $\mu^{(p)}= U([0,1])$, for every $p\!\in (0, +\infty)$. 

\smallskip
By {\em asymptotically $(L^p,\mu) $-optimal},  we mean that the $(L^p,\mu)$-mean quantization errors induced by the grids $\Gamma_{n}$ satisfy the sharp asymptotics of Zador's Theorem, namely $\lim_n N_n^{\frac 1d}  e_{p}(\Gamma_n,X) =  \widetilde J_{p,d} \|\varphi\|_{L^{\frac{d}{p+d}}(\lambda_d)}^{\frac{1}{p}}$.

\medskip
It is pointed out in~\cite{Brancoetal} (Theorem~4.10 and Corollary~4.11) that the quantizers $(a^{(N),p})_{N\ge 1}$ designed from an $(L^p,\mu)$-optimal greedy quantization sequence $(a_{_{N,p}})_{N\ge 1}$ are usually not asymptotically $(L^p,\mu)$-optimal, even up to an extraction.  The counter-example is exhibited in the $1$-dimensional basic  setting  of the uniform distribution $U([0,1])$.

The authors first build and analyze an  $(L^1, U([0,1]))$-optimal greedy sequence $(a_{_{N,1}})_{N\ge 1}$. Then, they show    that the (tight) sequence of empirical measures $\widetilde \mu_{_N} = \frac 1N\sum_{1\le k\le N} \delta_{a_{k,1}}$  on $[0,1]$ 
does not have the uniform distribution $U([0,1])$ (or equivalently the Lebesgue measure $\lambda_{1|[0,1]}$ over $[0,1]$)  as a weak limiting distribution. In particular, this implies, owing to  the above empirical measure theorem, that
\[
\liminf_N  N e_1\big(a^{(N),1},U([0,1])\big)>\frac 14 = \widetilde J_{1,1}
\]
keeping in mind that  $\displaystyle \widetilde J_{1,1}= \lim_N N e_{1,N}\big(U([0,1])\big)= \inf_N N e_{1,N}\big(U([0,1])\big)$. Otherwise, by the above empirical mean theorem, there would exist a subsequence $N'\to+ \infty$ such that $\widetilde \mu_{N'}\stackrel{(w)}{\rightarrow} \mu^{(1)}=\lambda_{1|[0,1]}= U([0,1])$. Equivalently, this reads 
\[
\liminf_N \frac{e_1\big(a^{(N),1},U([0,1])\big)}{e_{1,N}\big(U([0,1])\big)}>1.
\]
Numerical tests graphically reproduced in~\cite{Brancoetal} (Figure~1, p.521) suggest that  
$$
\displaystyle \liminf_N  Ne_{1}\big(a^{(N),1},U([0,1])\big)  \approx0.255 \approx 1.02 \times \widetilde J_{1,1}.
$$
 Our own numerical tests, based on the algorithms developed in Section~\ref{sec:greedyalgo} in the quadratic case ($p=2$), implemented   with  the uniform distribution,  the scalar ${\cal N}(0,1)$ and     bi-variate ${\cal N}(0;I_2)$ normal distributions  provide  similar conclusions  (see~Section~\ref{sec:greedyalgo} devoted to algorithmic aspects and numerical experiments).

\smallskip
This leads to our first open question: is this a generic situation? Or, to be more precise:

\medskip
\noindent   {\bf Open question~1}: May an optimal $(L^p,\mu)$-greedy sequence $(a_{_{N,p}})_{N\ge 1}$    contain  subsequence(s) $\big(a^{(N'),p}\big)_{N\ge 1}$ of asymptotically $(L^p,\mu)$-optimal $\mu$-quantizers?

\bigskip In fact, we conjecture that the a generic answer  is negative. This amounts to proving,  still owing to the empirical measure theorem,  that for any optimal $(L^p,\mu)$-greedy sequence $(a_{_{N,p}})_{N\ge 1}$
\[
\liminf_N \frac{e_p\big(a^{(N),p},\mu\big)}{e_{p,N}\big(\mu\big)} >1.
\]
\paragraph{$\rhd$ Are $(L^p,\mu)$-optimal greedy quantization sequence really optimal among ($\mu$-rate optimal) sequences?} Let us have a look at the celebrated dyadic Van der Corput ({\em VdC}$\,$) sequence, viewed as a quantization sequence.  Let us recall that the dyadic  {\em VdC}$\,$ sequence is defined by
\[
\xi_N = \sum_{k=0}^r \frac{n_k}{2^{k+1}}\quad \mbox{ where}\quad N=n_r2^r+\cdots+n_0,\; n_i \!\in \{0,1\},\, i=1,\ldots,r.
\] 

\smallskip
\noindent $\rhd$   {\em The $L^1$-mean quantization problem for the  {\em VdC}$\,$ sequence.} Elementary computations, not reproduced here, show  that
\[
\liminf_N N e_1\big(\xi_1,\ldots, \xi_{_N},[0,1]\big)= \frac{1}{4} = \widetilde J_{1,1}
\]
and that
\[
\limsup_{N}N e_1\big(\xi_1,\ldots, \xi_{_N},[0,1]\big)= \frac{9}{32}= \frac 98\, \widetilde J_{1,1}.
\]

This $\liminf$ is achieved by the subsequence $N_n= 2^{n-1}$, $n\ge 1$, and the $\limsup$ with subsequence $N_n= \frac 32. 2^n= 3. 2^{n-1}$, $n\ge 1$. So we can claim that: 

\begin{itemize}
\item there exist rate optimal sequences in the sense of~\eqref{eq:RateOpt} which are not solutions to the greedy problem~\eqref{eq:greedy};
\item there exist rate optimal sequences $(\xi_N)_{N\ge 1}$ containing subsequence of quantizers $(\xi^{(N')})_{N\ge 1}$ which are asymptotically $L^1$-rate optimal quantizers: so is the case of the {\em VdC}$\,$ sequence with the above subsequence $N'= 2^{n-1}$.
\end{itemize}

Figure~1 in~\cite{Brancoetal} also suggests that the $L^1$-optimal greedy quantization sequence $\big(a_{N,1}\big)_{N\ge 1}$ for the uniform distribution $U([0,1])$ satisfies 
$$
\displaystyle \limsup_N Ne_{1}\big(a^{(N),1}, U([0,1])\big ) \approx 1.09 \times \widetilde J_{1,1}\; \mbox{ and }\; 1.09 <1.125 = 9/8.
$$

\noindent 
$\rhd$ {\em The  $L^2$-mean quantization problem for the  {\em VdC}$\,$ sequence.}  The same phenomenons are confirmed in the quadratic case    since,  {\em mutatis mutandis},
\[
\liminf_N N e_2\big(\xi_1,\ldots, \xi_{_N},[0,1]\big)= \frac{1}{2\sqrt{3}} =\widetilde J_{2,1}\; \mbox{ and }\; 
\limsup_{N}N e_2\big(\xi_1,\ldots, \xi_{_N},[0,1]\big)=   \frac{3\sqrt{5}}{4}\times \widetilde J_{2,1}
\]
where we keep in mind that $\displaystyle  \widetilde J_{2,1} =\lim_N Ne_{2,N}(U([0,1])) =\inf_NNe_{2,N}(U([0,1]))$. 

\begin{figure}\label{fig:Unif10000}
  \centering
  
  \vskip -5 cm 
    \includegraphics[width=0.7\textwidth]{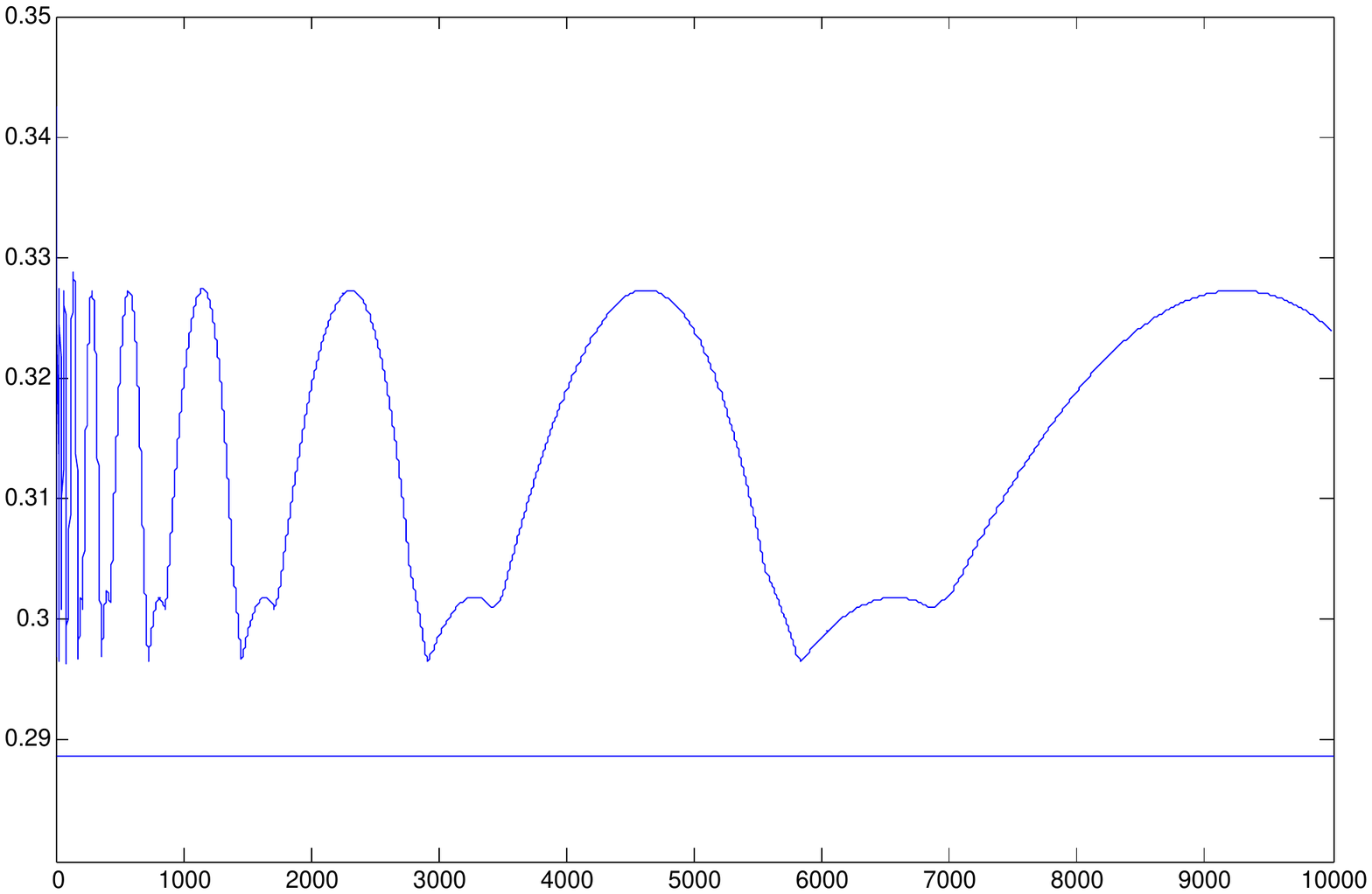}

     \vskip -5 cm 
      \caption{\it Graph $N\mapsto  Ne_{2}\big(a^{(N)}, U([0,1])\big)$, $N=1,\ldots, 10\,000$.}
  \end{figure}

\smallskip
On the other hand,  in a quadratic framework, using the greedy Lloyd~I procedure   described and analyzed  in the next Section~\ref{sec:greedyalgo1} (see Equations~\eqref{eq:greedyLloyd} if $d=1$ and~\eqref{eq:Lloyd-I-d} if $d\ge 2$), we   also observe numerically (see Figure~\ref{fig:Unif10000}) that 
\[
  \liminf_N Ne_{2}(a^{(N),2},U([0,1]) )\approx  0.29656\approx 1.02732  \times \widetilde J_{2,1}>\widetilde J_{2,1}
  \]
 and
  \[
   \limsup_N Ne_{2}\big(a^{(N),2},U([0,1]) \big) \approx 0.32736 \approx  1.13401  \times  \widetilde J_{2,1}
   \]
since $\widetilde J_{2,1} = \frac{1}{2\sqrt{3}}$. So the ``loss" is about $13\,\%$.

As for the $\liminf$, we verify  again  that no subsequence of $\big (a^{(N),2}\big)_{N\ge 1}$  can be asymptotically $L^2$-optimal and, ss for the $\limsup$,   that the quadratic optimal greedy sequence $\big(a_{N,2}\big)_{N\ge 1}$   outperforms the dyadic {\em VdC}$\,$ sequence from the $\limsup$ criterion  since $1.13401  <\frac{3\sqrt{5}}{4}=1.67706$. 

\medskip
\noindent $\rhd$ {\em Concatenated sequences.} From a more general point of view, there is a canonical method to produce for any distribution $\mu$ on $(\R^d, {\cal B}or(\R^d))$,  a $\mu$-rate optimal sequence for $(L^p,\mu)$-quantization by {\em concatenating}    $(L^p,\mu)$-optimal grids of size $2^\ell$. We proceed as  follows. 
%
Let   $(b_N)_{N\ge 1}$ be a sequence made up with $(L^p,\mu)$-optimal quantizers at level $2^{\ell}$, $\ell=0,\ldots n-1$ $i.e.$ so that
 \begin{equation}\label{eq:concatseq}
 \Big\{b_{2^\ell}, \ldots,b_{2^{\ell+1}-1}\Big\}\; \mbox{ is an $(L^p,\mu)$-optimal quantizer at level } 2^{\ell}.
 \end{equation}
 One checks straightforwardly by monotony of the $L^p$-mean quantization error that, for every $n\ge 1$, 
 \[
 e_{2^{n}-1}(b^{(2^{n}-1)},\mu) \le e_{2^{n-1}}\big(\{b_{2^{n-1}},\ldots,b_{2^{n}-1}\},\mu\big).
 \]
 Hence, for every $N\ge 1$, let be $n=n(N)$ be such that $2^{n}-1\le N\le 2^{n+1}$. Then 
 \[
 e_p(b^{(N)}, \mu)\le e_p(\{b_{2^{n-1}}, \ldots,b_{2^n-1}\})= e_{p, 2^{n-1}}(\mu)
 \]
 so that 
 \[
 \limsup_N N^{\frac 1d}e_p(b^{(N)}, \mu)\le \limsup_N \left(\frac{N}{2^{n(N)}}\right)^{\frac 1d} \lim_N N^{\frac 1d} e_{p,N}(\mu)= 2^{\frac 1d}  \lim_N N^{\frac 1d} e_{p,N}(\mu).
 \]
 
\noindent$\rhd$ {\em First elements of comparison.}

\smallskip 
-- If $\mu=U([0,1])$ and $ p=1$, one easily checks by induction that  the dyadic  {\em VdC}$\,$ sequence can be obtained as a properly reordered sequence $(b_N)_{N\ge 1}$  from the $L^p$-optimal quantizers at level $N$ given by $\Big\{\frac{2k-1}{2N},\, 1\le k\le N\Big\}$ when $N=2^n$, $n\ge 0$. In this very situation, the factor $2^{\frac 1d}=2$ is conservative since it can be replaced when $p=1$ by $\frac 98= 1.125$ as seen above. 

Anyway,  the  $L^1$-optimal greedy  quantization sequence  keeps the lead, since  $\displaystyle \limsup_N \frac{e_1(a^{(N),1}, \mu)}{e_{1,N}(\mu)}\approx 1.09< \frac{9}{8}\approx 1.125$.

\smallskip 
--   If $\mu=U([0,1])$ and $ p=2$, once again,  the  quadratic  optimal  greedy quantization sequence again keeps the lead, since  
$$
  \limsup_N \frac{e_2(a^{(N),2}, \mu)}{e_{2,N}(\mu)}\approx 1.13401< \frac{3\sqrt{5}}{4}\approx 1.67706<2.
  $$

\smallskip 
--  If $\mu = {\cal N}(0;I_2)$ (bivariate normal distribution   $i.e.$ $d=p=2$), our own numerical experiments suggest  for the third time   (see more detailed numerical results in Section~\ref{ssec:N0I2}) that a quadratic optimal greedy  quantization sequence (or, in practice,   the  suboptimal sequence resulting from the  numerical implementation of  the greedy Lloyd~I algorithm) has a lower constant than $2^{\frac 1d}\times   \lim_N N^{\frac 12} e_{2,N}{\cal N}(0 ; I_2)$.

\smallskip
All these considerations  experiments lead us to formulate a second  open question:

\medskip
\noindent {\bf  Open question~2}: Does an $(L^p,\mu)$-optimal greedy  quantization   $(a_{_{N,p}})_{N\ge1}$ produce the lowest value for  $\displaystyle \limsup_N N^{\frac 1d}\,e_{p,N}\big(a^{(N),p},\mu\big)$ among all sequences $(a_{_{N,p}})_{N\ge 1}$~?

\medskip A less ambitious question could be to compare $(L^p,\mu)$-optimal greedy sequences to  concatenated sequences~\eqref{eq:concatseq} $i.e.$: ``Is  the (strict) inequality $\displaystyle  \limsup_N \frac{e_{p}\big(a^{(N),p},\mu\big)}{e_{p,N}(\mu)}< 2^{\frac 1d}$ always satisfied?"

\paragraph{$\rhd$ Practical aspects in view of numerics.} From a more applied point of view, it would be of inte\-rest to establish  for $(L^p,\mu)$-optimal greedy sequences a counterpart of the non-asymptotic Zador Theorem in order to  upper-bound the $(L^p,\mu)$-mean quantization error of any greedy optimal sequence (normalized by $N^{-\frac 1d}$) by the $L^{p+\delta}$-pseudo-standard deviation of the distribution $\mu$ and a universal constant depending only on $d,\,p$ and $\delta$. The proof of the non-asymptotic Zador's Theorem (a slight improvement of  Pierce's Lemma established $e.g.$ in~\cite{LUPAaap}) relies on a random quantization argument involving the random quantizers $(Y^{(N)})_{N\ge 1}$ designed from an i.i.d. sequence $(Y_N)_{N\ge 1}$ with an appropriate distribution $\nu$, such a result is not hopeless.

\smallskip For numerical purposes, in particular numerical integration or conditional expectation approximation, some reasonably good estimates of $\limsup_N N^{\frac 1d}e_{1,N}(a^{(N)}, \mu)$ in~\eqref{eq:RateOpt} would be very useful. This is to be compared to the never ending quest for sequences with low discrepancy with  lower constant in the  Quasi-Monte Carlo community.   

 \section{Algorithmic aspects in the quadratic case}\label{sec:greedyalgo}   In this section we assume that $\R^d$ is equipped with the canonical Euclidean norm and that $p=2$ (purely quadratic setting).  So, will simply denote $(a_{_N})_{N\ge1}$ quadratic  optimal greedy sequences. 
 
 Practical computation of an optimal greedy sequence of quantizers relies on obvious variants  algorithms ($CLVQ$ and Lloyd) implemented recursively: to switch from $N$ to $N+1$, one first  adds a $(N+1)^{th}$ point (sampled from the support of the distribution $\mu$)  to the $N$-tuple  $(a_1,\ldots a_{_N})$ computed during the first  $N^{th}$ stages of the optimization  procedure. This makes  the starting $(N+1)$-tuple for the modified $CLVQ$ to Lloyd procedure. Then, one launches  one of these two optimization procedures  with the following restriction:   {\em all formerly computed components $a_i$, $1\le i\le N-1$ are kept frozen}, and only the new point is moved following the standard rules. Thus, if implementing a  $CLVQ$ like procedure, when the $N^{th}$ component is the ``winner" in the competition phase ($i.e.$ the $N^{th}$ component is the nearest neighbour to the  new  input stimulus). As for the (randomized) Lloyd~I procedure,    the Voronoi cell of the $N^{th}$ component  is the only one whose centroid (the $N^{th}$ component)  is updated, the other $N-1$ components remaining frozen as well.  Let us be   more precise. 
 
  \subsection{The one-dimensional quadratic case}\label{sec:greedyalgo1} When $d=1$  and  the distribution $\mu$ is absolutely continuous with a continuous positive probability density $\varphi$ on the real line, one can directly consider  the counterpart of the historical deterministic Lloyd~I procedure and of the gradient descent sometimes known as  Forgy's algorithm or $k$-means. Let us be more specific.
 
 \paragraph{$\rhd$ Greedy Lloyd's~I procedure}
 
 $\bullet$ Assume $a_1,\ldots a_{N-1}$ have been computed. Let $a^{(N-1)}_1<\cdots < a^{(N-1)}_{N-1}$ be  an  increasing reordering of $a_1,\ldots,a_{_{N-1}}$. 
 
 \smallskip
 \noindent $\bullet$  Assume   the $N$  inter-point  local inertia has also been computed, namely
  \[
\sigma^2_i:= \int_{a^{(N-1)}_{i}}^{a^{(N-1)}_{i+\frac 12}}|a^{(N-1)}_i-\xi|^2\mu(d\xi)+\int_{a^{(N-1)}_{i+\frac 12}}^{a^{(N-1)}_{i+1}}|a^{(N-1)}_{i+1}-\xi|^2\mu(d\xi),\, \; i=1,\ldots,N
 \]
 where 
 \[
 a^{(N-1)}_{0}=  a^{(N-1)}_{\frac 12} = -\infty,\; a^{(N-1)}_{i-\frac 12}= \frac{a^{(N-1)}_{i-1} + a^{(N-1)}_i}{2},\ i=2,\ldots,N-2, \, a^{(N-1)}_{N-\frac 12}= a^{(N-1)}_{N}= +\infty.
 \]

 \smallskip
 \noindent $\bullet$ Choose an index $i_0=i_0(N-1)$ such that $\sigma_{i_0}^2 =\max_{0\le i\le N} \sigma^2_i$ (maximal local inertia), then  consider $a_0= a_{N,0}\!\in (a^{(N-1)}_{i_0},a^{(N-1)}_{i_0+1})$ and  finally define recursively  a sequence $a_{[n]}= a_{N,n}$, $n\ge 1$, by
 \begin{equation}\label{eq:greedyLloyd}
 a_{[n+1]}=\E\big(X\,|\, X\in W_{N,[n]}\big)= \frac{K_{\mu}\big(\frac{a^{(N-1)}_{i_0+1}+a_{[n]}}{2}\big)-K_{\mu}\big(\frac{a^{(N-1)}_{i_0}+a_{[n]}}{2}\big)}{F_{\mu}\big(\frac{a^{(N-1)}_{i_0+1}+a_{[n]}}{2}\big)-F_{\mu}\big(\frac{a^{(N-1)}_{i_0}+a_{[n]}}{2}\big)},\; n\ge 0,
 \end{equation}
 where $F_{\mu}(x)= \mu((-\infty,x])$ is the cumulative distribution function of $\mu$ and $K_{\mu}$ its cumulative first moment function defined by 
 \[
 K_{\mu}(x)=\int_{(-\infty,x]}\!\! \xi\,\mu(d\xi),\;\; x\!\in \R.
 \]
 
It follows form an easy induction that, at  every step $n\ge 0$ of the  procedure, $a_{N,[n]}\!\in W_{N,[n]}\subset  (a^{(N-1)}_{i_0},a^{(N-1)}_{i_0+1})$ so that the procedure is well-defined.

\begin{Pro}\label{pro:Greedy-Lloyd} If $\mu$ is {\em strongly unimodal } in the sense that $\mu= \varphi.\lambda_1$ with $\varphi:\R\to\R$  $\log$-concave, then  $a_{N,[n]}$ converges toward the unique solution $a_{N,\infty}\!\in (a^{(N-1)}_{i_0},a^{(N-1)}_{i_0+1})$ of the fixed point equation  
\begin{equation}\label{eq:Stationarity}
a_{_N}= \E\big(X\,|\, X\in W_{N}\big)
\end{equation}
where $W_{N}\subset  (a^{(N-1)}_{i_0},a^{(N-1)}_{i_0+1})$  is the closed Voronoi cell of $a_{_N}$ in $a^{(N-1)}\cup\{a_{_N}\}$.
\end{Pro}
 
The detailed proof is postponed to the Appendix~\ref{app:A1}. But we can already mention that it relies on classical arguments called upon  in the proofs of the convergence of the standard Lloyd~I procedure (and the uniqueness of the possible stationary limiting point, see~\cite{KIE, BOPA1}).

 \medskip
\noindent {\bf Remarks.} $\bullet$  The computation of the integrals involved in the algorithm can be performed by higher order quadrature formulas, or $e.g.$ in the case where  $\mu={\cal N}(0;1)$ using the closed form for $\int_{-\infty}^x \xi e^{-\frac{\xi^2}{2}} \frac{d\xi}{\sqrt{2\pi}}= -\frac{e^{-\frac{x^2}{2}}}{\sqrt{2\pi}}$ and high accuracy approximations for its cumulative distribution function  $\Phi_0$, using $e.g.$ continuous fractions expansions (see~\cite{Handbook}).
 
 \smallskip
 \noindent $\bullet$ The $\log$-concave assumption which implies the uniqueness of the fixed point for Equation~\eqref{eq:greedyLloyd}, is satisfied by many usual families of distributions on the real line like $e.g.$ the normal distributions ${\cal N}(m;\sigma^2)$, the exponential and Laplace distributions, the $\gamma( \alpha,\beta)$-distributions, $\alpha\ge 1$, $\beta>0$, are strongly unimodal. On the other hand, the Pareto distributions are not strongly unimodal though uniqueness holds true (see~\cite{FOPA}).
 %
 %
 
 \paragraph{$\rhd$ Greedy Forgy's algorithm (Newton zero search algorithm)} This procedure is  defined recursively by 
 \begin{equation}\label{eq:greedyForgy}
 a_{[n+1]}  =  a_{ [n]} - \Big(\gamma_{n+1}\wedge\frac{1}{ \rho(a_{[n]}) }\Big) \int_{\frac{a^{(N-1)}_{i_0}+a_{[n]}}{2}}^{\frac{a^{(N-1)}_{i_0+1}+a_{[n]}}{2}}\big( a_{[n]}-\xi\big)\mu(d\xi)
 \end{equation}
 where $\gamma_{n+1}\!\in (0,1)$ goes to $0$ as $n\to+\infty$, $\sum_n\gamma_n=+\infty$ and  
 \[
 \rho(a)= \mu\Big(\Big[\frac{a^{(N-1)}_{i_0}+a}{2}, \frac{a^{(N-1)}_{i_0+1}+a}{2}\Big]\Big)+\frac{a-a^{(N-1)}_{i_0}}{2}f\Big(\frac{a+a^{(N-1)}_{i_0}}{2}\Big) +\frac{a^{(N-1)}_{i_0+1}-a}{2}f\Big(\frac{a+a^{(N-1)}_{i_0+1}}{2}\Big)>0 
 \]
 is the second derivative of the a function $a\mapsto \E\big(\min|X-a_i|^2\wedge |X-a|^2\big)$.

 \smallskip
Note that, owing to the thresholding of $1/\rho(a_{N,[n]})$ by $\gamma_{n+1}\!\in (0,1)$, this procedure lives in the interval $(a^{(N-1)}_{i_0}, a^{(N-1)}_{i_0+1})$ which makes it  well-defined and consistent for every $n$.  
%
%

 \smallskip
When $\mu$ is not absolutely continuous,  one can implement the same procedure by removing the term involving the second derivative  with a step $\gamma_n$ satisfying the standard {\em decreasing step assumption} ($\sum_n \gamma_n=+\infty$ and $\sum_n \gamma^2_n<+\infty$), provided one can compute the $\mu$-integrals of interest.

\paragraph{$\rhd$ Numerical illustration with the ${\cal N}(0;1)$ distribution} To compute a quadratic  optimal greedy sequence of the normal distribution  $\mu = {\cal N}(0;1)$, we will take advantage of its symmetry. To this end we consider the distribution $\tilde \mu= \mu(\,.\,|\; \R_+)$ ($\mu$ conditioned to stay non-negative) which is clearly strongly unimodal and we compute by induction its quadratic optimal greedy sequence $(\widetilde a_{_N})_{N\ge 1}$ by the greedy Lloyd~I procedure~\eqref{eq:greedyLloyd}   with the convention that the origin $0$ is a fixed but {\em active} point as a possible nearest neighbour  for this slight variant. To be precise, we mean that $0$ has its own Voronoi cell in $\R_+$ or, equivalently, that we implement the algorithm,  starting at $\widetilde a_0=0$ when  $N=0$. 

As a second step, it is straightforward that the sequence defined by 
\[
a_0=0,\; a_{_{2N-1}}= \widetilde a_{_N},\; a_{_{2N}}= - \widetilde a_{_N},\; N\ge 1,
\]
is a quadratic optimal  greedy  sequence.

\begin{figure}\label{fig:Gauss10000}
  \centering
  
  \vskip -5 cm 
    \includegraphics[width=0.75\textwidth]{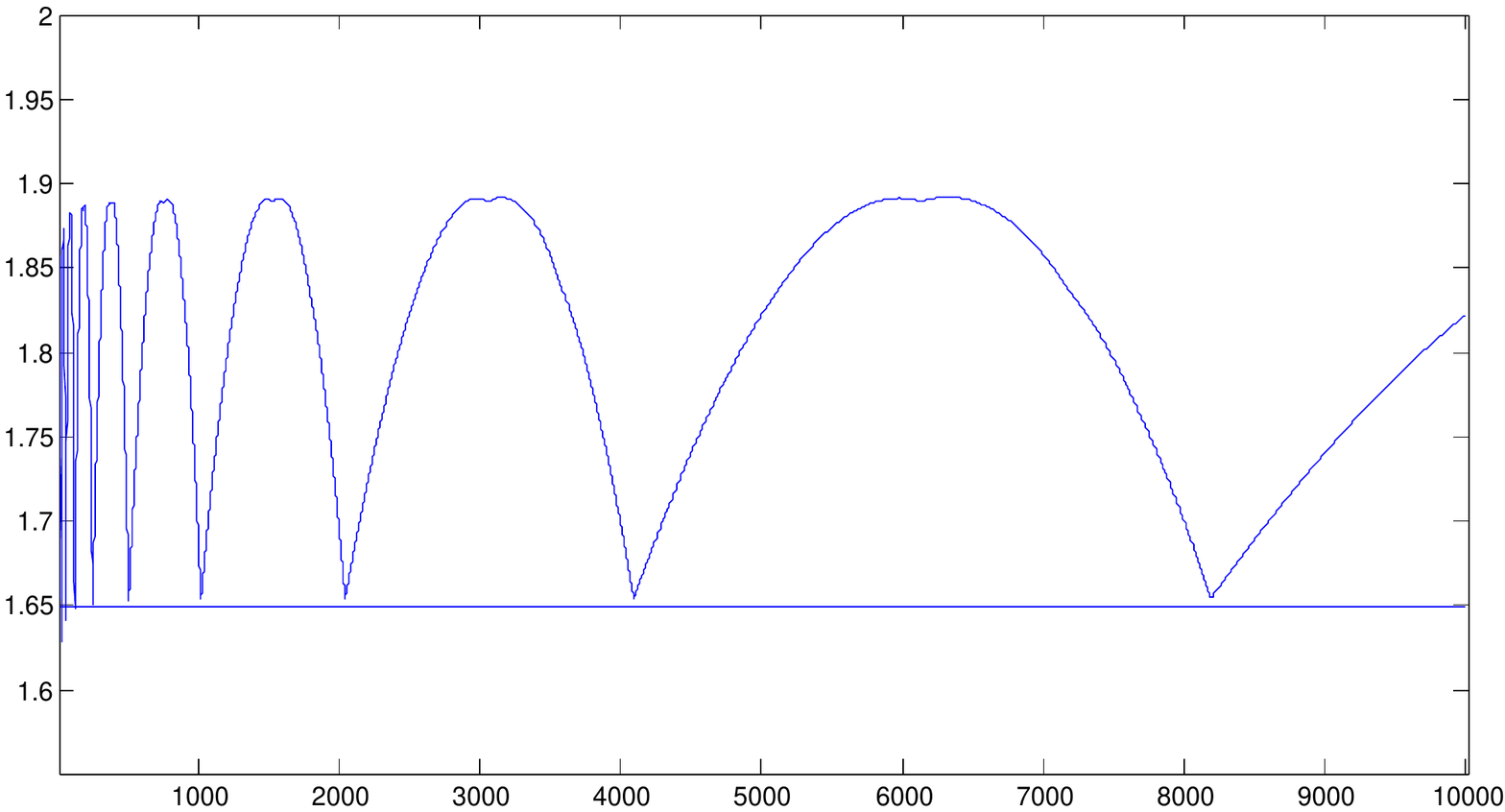}

     \vskip -5.5 cm 
      \caption{\it Graph $N\mapsto (2N-1)e_{2}\big(a^{(2N-1)}, {\cal N}(0;1)\big)$, $N=4,\ldots, 10\,000$.}
  \end{figure}

\smallskip We reproduce in Figure~2
the graph $N\mapsto (2N-1)e_{2}\big(a^{(2N-1)}, \mu\big)$, $N=4,\ldots,2^{10}=10\,000$, where $\mu={\cal N}(0;1)$.
 
  Note that $\displaystyle  \limsup_N Ne_{2}\big( a^{(N)}, \mu\big) = \limsup_N (2N-1)\, e_{2}\big(a^{(2N-1)},\mu\big)$ since $e_{2}\big( a^{(N)},\mu\big) \downarrow 0$ as $N\to+\infty$. As a consequence, we derive that 
 \[
\liminf_N Ne_{2}\big( a^{(N)}, {\cal N}(0;1)\big)\approx 1.6534\dots>  \sqrt{\frac 32 }\,\pi^{\frac 14}  = \lim_N Ne_{2,N}\big( {\cal N}(0;1)\big)
 \]
since $  \sqrt{\frac 32 }\,\pi^{\frac 14} \approx  1.63055 $. (The real constant in the right hand side of the inequality easily follows from Zador's Theorem). Note that, for the values   $N=2^n$, $0\le n\le 7$, we observe that $Ne_{2}\big( a^{(N)}, {\cal N}(0;1)\big)< \sqrt{\frac 32 }\,\pi^{\frac 14} $~(\footnote{This is consistent in some way with the conjecture that $Ne_{2,N}\big( {\cal N}(0;1)\big)$ is increasing toward its limit.}).

\smallskip
As for the limsup, we observe numerically that
 \[
\limsup_N Ne_{2}\big( a^{(N)}, {\cal N}(0;1)\big)\approx  1.8921<2\times  \sqrt{\frac 32 }\,\pi^{\frac 14}\approx 3.2611.
 \]
 Consequently, the highest ``loss" for this one-dimensional distribution with unbounded support is approximately of $15.7\,\%$. 
\subsection{The multidimensional quadratic case (higher dimensions)}\label{ssec:N0I2}
 
 In higher dimensions, deterministic procedures like deterministic greedy Lloyd's~I (fixed point procedure defined by~\eqref{eq:greedyLloyd})  or the greedy Forgy's (recursive zero search defined by~\eqref{eq:greedyForgy})  algorithms become computationally too demanding due to the repeated computations of integrals on the Voronoi cells of the quantizers. So, it becomes necessary, at least when $d\ge 3$,  to switch to stochastic optimization procedures like those described  below, which are adaptations of the stochastic procedures introduced to compute true optimal $N$-quantizers. For more details about  these original stochastic optimization procedures, mostly devised   in the 1950's, we refer $e.g.$ to~\cite{BMP, PAPR} for $CLVQ$ and~\cite{KIE, Du, PAYU} for (randomized) Lloyd's~I procedure or more applied textbooks like~\cite{GEGR}. These procedures have been extensively implemented to compute for numerical probability purposes  optimal grids of  $d$-dimensional normal distributions ${\cal N}(0,I_d)$ for  $d=1,\ldots,10$ and sizes up to $N= 10\,000$. 
 
 From a theoretical point of view, the common feature of these stochastic algorithms is that the convergence results ($a.s.$ or in $L^p$) remain partial, especially little is known when the distribution $\mu$ is not compactly supported.  So we present below their greedy variants (without rigorous proof as concerns   $CLVQ$). From a practical point of view,  for both procedures, the computation of  integrals on the Voronoi cells is replaced by repeated nearest neighbor searches among the components of the current $N$-quantizers which make them rather slow. But in our greedy framework, this drawback could be overcome by appropriate localization around the elementary quantizer of interest. But this is beyond the scope of the present work. 

%
%
\bigskip 
\noindent $\rhd$ {\sc (Randomized) greedy Lloyd's~I like procedure.} The {\em greedy Lloyd~I procedure} to compute $a_{_N}$, assuming that $a^{(N-1)}$ is known, (starting from the mean $a_1= \E \,X$)  can be recursively defined in the quadratic case as follows:
 \begin{equation}\label{eq:Lloyd-I-d}
  a_{N,[n+1]}= \E\big(X\,|\, X\!\in W_{N,[n]}\big),\quad a_{N,[0]} \!\in \R^d\!\setminus\{a^{(N-1)}\},
 \end{equation}
 where $W_{N,[n]}$ is the closed Voronoi cell of $a_{N,[n]}$ with respect to the quantizer $a^{(N-1)}\cup\{a_{N,[n]}\}$. Of course in practice, we stop the Monte Carlo simulation at finite range $M_n$. 
 
\bigskip
We establish in the proposition below, at least for absolutely continuous distributions with  convex support, that 
 \[
 \lim_{n\to+\infty} a_{N,[n]} \;\mbox{ does exists} 
 \]
under a local finiteness assumption on the possible equilibrium points.   Due to the existence of several equilibrium points, especially  in higher dimension, this limit may  not be  the solution to the greedy optimization problem at level $N$, but only a local minimizer. However, in practice, it turns out to be a good candidate.  

 
 \begin{Pro}\label{pro:LloydI-d} Assume the distribution $\mu$ of $X$ is strongly continuous ($i.e.$ assigns no mass to hyperplanes) with a convex support denoted $C_{\mu}= {\rm supp}(\mu)$. Then the above sequence $(a_{N,[n]})_{n\ge 0}$ is bounded and there exists  $\ell\!\in\big [e_{2}(a^{(N)}), e_{2}(a^{(N-1)}\cup\{a_{[0]}\}\big)$ such that  the set ${\cal A}_{\infty}(a_{[0]})$ of its limiting points   is a connected compact  subset of  the set $\Lambda_{\ell}$ of $\ell$-{\em stationary points}  defined by
 $$
 \Lambda_{\ell} = \Big\{a\!\in \R^d\,|\, e_{2,N}\big(a^{(N-1)}\cup\{a\}\big)=\ell \; \mbox{ and }\; a= \E\big(X\,|\, X\!\in W_{N,a}\big)\Big\}
 $$
 where $W_{N,a}$ denotes the closed Voronoi cell of $a$ induced by the $N$-quantizer $a^{(N-1)}\cup\{a\}$. In particular, $e_{2}\big(a^{(N-1)}\cup\{a_{[n]}\},X\big)\to \ell$ as $n\to +\infty$.
 
Furthermore, if  the  $\ell$-{\em stationary set} $\Lambda_{\ell}$ is locally finite ($i.e.$ with a  finite trace on  compact sets of $\R^d$),  then $a_{N,[n]}$ $a.s.$ converges to some point in $\Lambda_{\ell}$.
 \end{Pro}
 
%
%
%
 The proof is postponed to  Appendix~\ref{app:A2}.
 
 \medskip
 The true algorithm  to be implemented in practice is a {\em randomized}  version of this procedure where each conditional expectation is computed by Monte Carlo simulation (provided $X$ can be simulated at a reasonable cost): let $(X^m)_{m\ge 1}$ be an i.i.d. sequence of copies of $X$ (with distribution $\mu$) defined on a probability space $(\Omega,{\cal A}, \P)$. Then, by the Strong Law of Large Numbers, 
 \[
 a_{N,[n+1]} =\lim_{M\to +\infty} \frac{\sum_{m=1}^M X^m\mbox{\bf 1}_{\{X^m\in W_{N,[n]}\}} }{\sum_{m=1}^M \mbox{\bf 1}_{\{X^m\in W_{N,[n]}\}}  }\quad \P\mbox{-}a.s.
 \]
 
 \medskip
 \noindent  $\rhd$  {\sc Sequential Competitive Learning Vector Quantization   procedure}: Let $(\gamma_n)_{n\ge 1}$ be a sequence of $(0,1)$-valued step parameters satisfying a so-called decreasing step assumption: $\displaystyle \sum_n \gamma_n=+\infty$ and $\sum_n \gamma^2_n<+\infty$. Then set 
 \[
 a_{N,[n+1]}  =  a_{N, [n]} -\gamma_{n+1} \mbox{\bf 1}_{\{|X^{n+1}-a_{N, [n]}|< \min_{a\in a^{(N-1)}}  |X^{n+1}-a|\}} \big( a_{N,[n]} -X^{n+1}\big),\; a_{N,[0]} \!\in \R^d.
 \]
One may conjecture and experimentally check, at least for distribution with compact convex support, 
 \[
 \lim_{n\to+\infty} a_{N,[n]}= a_{_N}.
 \]
 If so is the case, one may apply the so-called {\em Ruppert-Polyak principle} which states that choosing a ``slowly decreasing" step of the form $\gamma_n =\frac{c}{c+n^{\alpha}}$, $\frac 12< \alpha<1$, and averaging the procedure by setting
 \[
 \bar a_{N,[n]}= \frac 1n \big(a_{[N,0]}+\dots+a_{[N,n-1]}\big), \; n\ge 1,
 \]
 will speed up the convergence or, to be more precise, will satisfy a Central Limit Theorem at rate $\sqrt{n}$ with the lowest possible asymptotic variance (see~$e.g.$~\cite{LUS, PAGSpring} for  details).
 
\medskip
 \noindent $\rhd$  {\sc Randomized Greedy Lloyd's~I randomized procedure for the bi-variate normal distribution} 

Let $\mu ={\cal N}(0;I_2)$ be the bi-variate normal distribution on the plane. Figure~3
depicts the graph of $N\mapsto \sqrt{N} e_{2}\big(a^{(N)},\mu\big)$ for $N=1$ up to $1\, 000$ (and Figure~4
depicts $a^{(1000)}$). This  suggests  that this sequence remains  bounded.  However, we are not sure with such a rough procedure that the computed sequence $(a_N)_{N\ge 1}$ is the optimal greedy one: at each step/level,
 there are clearly many local parasitic minima and one should add, prior to computing  $a_N$, a pre-processing phase, like in 
 one dimension, in order to choose among  the areas  defined by  the Delaunay triangulation attached to $a^{(N-1)}$, the one which induces the minimal inertia. But this phase is numerically demanding and has not yet been included in the existing script. 
 
 The randomized greedy Lloyd's method~1  has been  implemented at each  level $N$ with $M= M(N)=1\,000\times N$ $i.i.d.$  simulations of the ${\cal N}(0;I_2)$ distribution. Owing to  Zador's Theorem, we know that optimal quadratic quantizers satisfy (asymptotically)   
 \[
 \lim_N \sqrt{N} e_{2,N}\big({\cal N}(0;I_2)\big)= 2\sqrt{2\pi}\widetilde   J_{2,2} = \frac 23\sqrt{\frac{5\,\pi}{\sqrt{3}}}\approx 2.0077
 \]
since, owing to~\cite{GRLU} (Theorem~8.15, p.120, and Examples 8.12, p.116, devoted to hexagon lattices), $\widetilde J_{2,2}=\frac 13 \sqrt{\frac{5}{2\sqrt{3}}}$. Consequently, the ``loss" is less than $10\,\%$. We verify on our own numerical experiments carried out with $N=1000$ that it is likely that 
$$
\sup_{1\le N\le 1000} \sqrt{N}\,e_2\big(a^{(N)},\,{\cal N}(0;I_2)\big) \precsim 
2.18.
$$ 
 As already mentioned,   it suggests again that the greedy quantization sequence outperforms the concatenated sequence~\eqref{eq:concatseq}  since  $2.18< \sqrt{2} \times  \frac 23\sqrt{\frac{5\,\pi}{\sqrt{3}}}\approx 2.8392 $ (even if one may guess that the factor $2^{\frac 1d}=\sqrt{2}$ is probably too conservative in practice).
 \begin{figure}\label{fig:NErreurQ}
  \centering
  
  \vskip -5cm
    \includegraphics[width=0.75\textwidth]{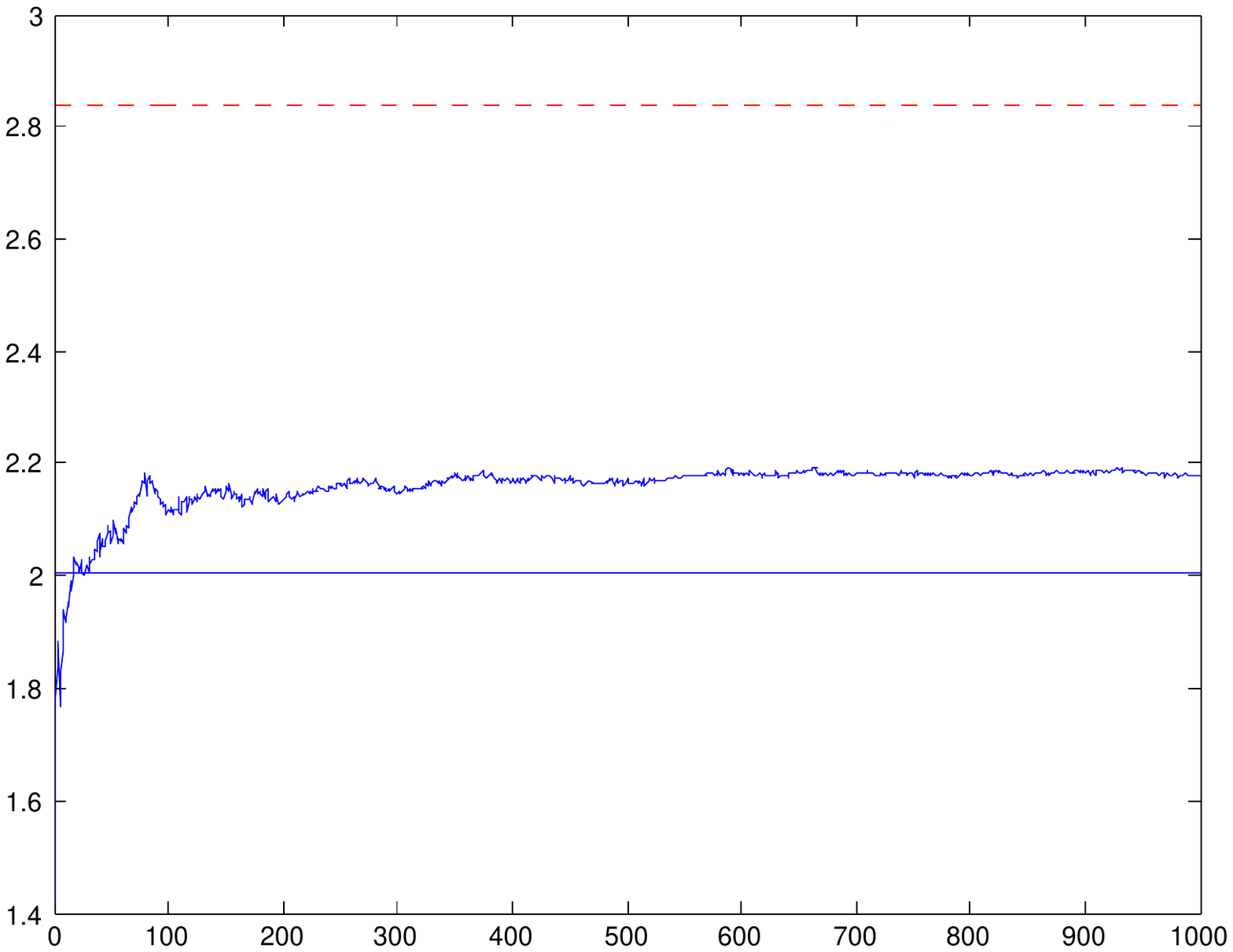}

      \vskip -5.25cm
      \caption{\it Graph $N\mapsto \sqrt{N}\,e_{2}\big(a^{(N)}, {\cal N}(0;I_2)\big)$, $N=1,\ldots,10^3$,    computed by the randomized greedy Lloyd~I procedure ($M=M(N)= 1\,000\times N$, $N=1,\ldots,10^3$). Flat solid line  ($\textcolor{blue}{-\!\!\!-\!\!\!-}$) depicts Zador's constant~$\widetilde   J_{2,2} = \frac 23\sqrt{\frac{5\,\pi}{\sqrt{3}}}$; flat dashed line ($- - -$) depicts the natural upper bound for the concatenated sequence.}
  \end{figure}

    \begin{figure}\label{fig:cloud1000}
  \centering

  \vskip -5cm
    \includegraphics[width=0.75\textwidth]{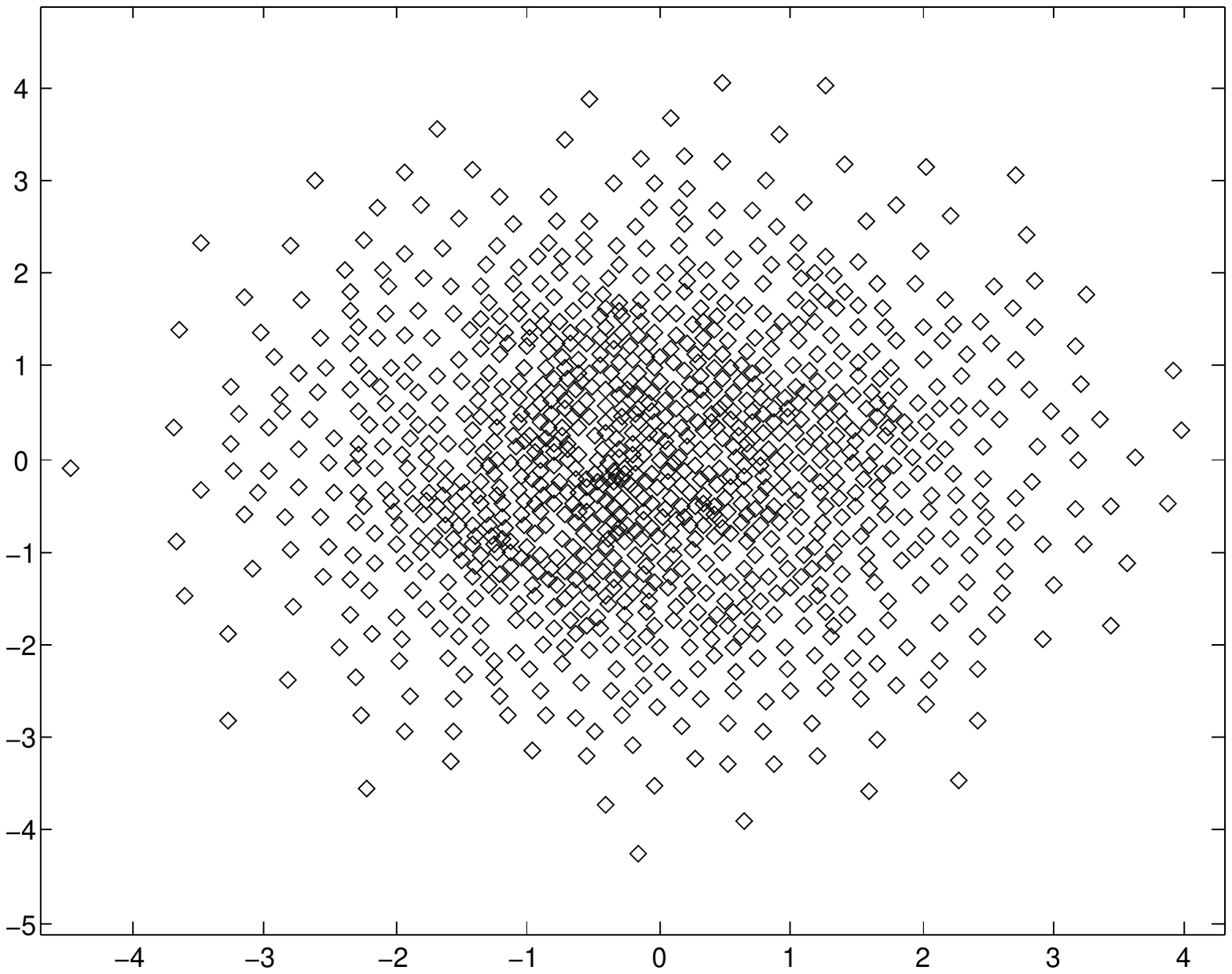}

  \vskip -5.25cm \caption{\it Greedy quantizer $a^{(1000)}$ for the ${\cal N}(0;I_2)$ distribution    computed by the randomized greedy Lloyd~I procedure with a simulation of size $M= 10^6$.}
  \end{figure}
  
\section{Greedy quantization versus Quasi-Monte Carlo?}\label{sec:GreedyvsQMC}
Of course, for every integer $N\ge1$,  the {\em weights}   induced by the $\mu$-mass of the Voronoi cells associated to $a^{(N)}$ define canonically a sequence of  $N$-tuples which usually cannot be ``arranged" into a sequence, even up to a re-scaling.   When considering the unit hypercube $[0,1]^d$ as a state space in $d$ dimension, it is easy natural to compare an optimal greedy sequence with respect to the uniform distribution $U([0,1]^d)$ and the so-called {\em uniformly distributed sequences} usually implemented in the Quasi-Monte Carlo method. 
 
 Let us recall that a  sequence $(\xi_{_N})_{N\ge 1}$ is uniformly distributed over $[0,1]^d$ if the empirical measures $\displaystyle \nu_{_N}=\frac 1N \sum_{i=1}^N \delta_{\xi_i}$ weakly converges toward the {\em Lebesgue measure} $\lambda_d$ on $[0,1]^d$. In particular this means that for every bounded $\lambda_d$-$a.s.$ continuous function $f:[0,1]^d \to \R$, $\displaystyle \frac 1N \sum_{i=1}^N f(\xi_i) \to \int_{[0,1]^d} \!\!fd\lambda_d =\int_{[0,1]^d}\!\! f(u)du$. This means that the weights associated to a uniformly distributed sequence  are by definition all equal to $\frac1N$ which leads to a simple   normalization factor $1/N$. What is  the cost induced by these uniform weights $\frac 1N$, compared to the optimal weights deduced from the cell (hyper-)volumes of the Voronoi diagram of $\xi_1,\ldots,\xi_{_N}$?  The answer is essentially $\log N$ and is provided by   Proinov's  theorem (see~\cite{PRO}) recalled below which evaluates precisely the convergence rate of empirical measures of uniformly distributed sequences on Lipschitz continuous  functions.

In the Quasi-Monte Carlo  ($QMC$) method, the performance of an $N$-tuple $(\xi_1,\ldots,\xi_{_N})\!\in ([0,1]^d)^N$ is measured by the Kolmogorov-Smirnov distance between the extended cumulative distribution function of its empirical measure $\nu_{_N},\,N\ge 1$, and the uniform distribution $U([0,1]^d)$, namely the so-called {\em star discrepancy} defined by
\begin{equation}\label{eq:Discstar}
D^*_N(\xi_1\ldots,\xi_{_N})=\sup_{u\in [0,1]^d}\left|\frac 1N \sum_{i=1}^N \mbox{\bf 1}_{\{\xi_i\in [\![0,u]\!]\}}- \lambda_d\big([\![0,u]\!]\big)\right|
\end{equation}
where $[\![0,u]\!]= \prod_{\ell=1}^d [0, u^\ell]$, $u= (u^1,\ldots,u^d)$.

\smallskip Several sequences $\xi= (\xi_{_N})_{N\ge 1}$ have been exhibited (see~\cite{NIE}) whose star discrepancy at the origin satisfies for  a real constant $C(\xi)\!\in (0,+\infty)$,
\begin{equation}\label{eq:discstar}
\forall\, N\ge 1,\quad D^*_N(\xi_1\ldots,\xi_{_N})\le C(\xi) \frac{(1+\log N)^d}{N}.
\end{equation}
Among them one can cite the $p$-adic $\!${\em VdC}$\,(p)\,$ sequences ($p\geq 2$ in $1$-dimension) and, when $d\ge 2$, the Halton sequences (whose $i^{th}$ component is the  $\!${\em VdC}$\,(p_i)\,$  sequence  where the bases $p_i$, $i=1,\ldots,d$, are the first $d$ prime numbers), the Faure sequences, the Sobol' sequences (a unifying framework has been developed by Niederreiter, see $e.g.$~\cite{NIE}). For definitions of these sequences and numerical tests on various problems we refer to~\cite{BOLE, PAXI}. Although such a rate has never been proved to be  the lowest possible, this opinion is commonly shared  by the $QMC$ community (however see again~\cite{NIE} or~\cite{PAGSpring} for a review of existing lower bounds).

\smallskip
The striking fact with these sequences satisfying~\eqref{eq:discstar},  called {\em sequences with low discrepancy}, is that when they are implemented on the class of functions with {\em  finite variation} on $[0,1]^d$ the  Koksma-Hlawka inequality implies that, for every such function $f:[0,1]^d \to \R$ 
\begin{equation}\label{eq:KokHla}
\left| \int_{[0,1]^d}\!\! f(u)du- \frac 1N \sum_{i=1}^N f(\xi_i)\right|\le V(f)D^*_N(\xi_1\ldots,\xi_{_N})
\end{equation}
where $V(f)$ denotes the variation of the function $f$. So it induces  for this specific class of functions a rate of numerical integration of order  $O\Big(\frac{(\log N)^d}{N}\Big)$. In one dimension ($d=1$), However,  the above notion of finite variation coincides with  the standard definition of finite variation in real analysis. 

When  $d\ge 2$, several definitions can be given, the most popular being the finite variation in the Hardy \& Krause sense (as described $e.g.$ in~\cite{NIE}).  Another slightly less general --~but more elementary~-- being the finite variation in the signed measure sense developed in~\cite{BOLE} (see also~\cite{PAGSpring}). Unfortunately, as the dimension $d$ increases, the set of functions with finite variation (in any of the above  senses) 
becomes somewhat ``sparse" among the set of all real-valued Borel functions defined on $[0,1]^d$. So this striking behavior may be considered as not significant when dealing with practical simulation problems. However to carry out a comparison, we need to evaluate their performances  the same significant functional space, namely  that of Lipschitz continuous functions.  Proinov's theorem below provides an answer. 

%

\begin{Thm}[Proinov~\cite{PRO}]\label{Proinov} Assume $\R^d$ is equipped with the $\ell^\infty$-norm $|(\xi^1,\ldots,\xi^d)|_{_\infty} = \max_{1\le i\le d}|\xi^i|$. For every continuous function $f:[0,1]^d\to \R$, we define uniform continuity modulus of $f$ (with range $\delta\!\in [0,1]$) by 
$$
w(f,\delta):= \sup_{\xi,\,\xi'\in [0,1]^d, \,|\xi-\xi'|_{_\infty} \le \delta}|f(\xi)-f(\xi')|.
$$ 
 
\ss \noindent $(a)$ Let
$(\xi_1,\ldots,\xi_{_N})\!\in ([0,1]^d)^N$. For every  continuous function $f:[0,1]^d\to \R$,
\[
\left | \int_{[0,1]^d} f(u)du-\frac 1N \sum_{i=1}^N f(\xi_i)\right|\le C_d\,w\big(f,D^*_{_N}(\xi_1,\ldots,\xi_{_N})^{\frac 1d}\big)
\]
where  $C_d\!\in(0,\infty) $ is a universal optimal real constant only depending on $d$. In particular, if the function $f:[0,1]^d \to \R$ is $\ell^{\infty}$-Lipschitz continuous with coefficient $[f]_{\rm Lip} := \sup_{x,y\in [0,1]^d}\frac{|f(x)-f(y)|}{|x-y|_{_\infty}}$, then
\[
\left|\int_{[0,1]^d} f(u)du-\frac 1N \sum_{i=1}^N f(\xi_ik\right|\le C_d\,[f]_{\rm Lip}D^*_{_N}(\xi_1,\ldots,\xi_{_N})^{\frac 1d}.
\]

If $d=1$, $C_d=1$ and if $d\ge 2$, $C_d\!\in [1,4]$.

\smallskip 
\noindent  $(b)$ In particular if $(\xi_N)_{N\ge 1}$ is a sequence with low discrepancy in the above sense, then
\[
\left|\int_{[0,1]^d} f(u)du-\frac 1N \sum_{i=1}^N f(\xi_)\right|\le C_d\,[f]_{\rm Lip}C(\xi)\frac{1+\log N}{N^{\frac 1d}}.
\]
\end{Thm}

\begin{Cor} \label{cor:KL-VQ}$(a)$ For every $N$-tuple $(\xi_1,\ldots,\xi_{_N})\!\in( [0,1]^d)^N$
\[
e_1\big(\xi_1,\ldots,\xi_{_N},U([0,1]^d)\big)\le C_dD^*_{_N}(\xi_1,\ldots,\xi_{_N})^{\frac 1d}.
\]

\noindent $(b)$ In particular, when $d=1$, $\displaystyle 
e_1\big(\xi_1,\ldots,\xi_{_N},U([0,1])\big)\le  D^*_{_N}(\xi_1,\ldots,\xi_{_N})$.
\end{Cor}
\noindent {\bf Proof (of $(b)$)}. Assume $d=1$. The function $f_{\xi}:u\mapsto \min_{1\le i\le N}|u-\xi_i|$ defined on $[0,1]$ is $1$-Lipschitz continuous,  hence has finite variation with $V(f_{\xi})=1$. Then Koksma-Hlawka Inequality~\eqref{eq:KokHla} or Proinov's error bound in $(a)$  both  imply that
 \begin{eqnarray*}
e_1\big(\xi_1,\ldots,\xi_N, U([0,1])\big)&= & \left| \frac 1N \sum_{i=1}^N f_{\xi}(\xi_i) -e_1\big(\xi_1,\ldots,\xi_N, U([0,1])\big)\right|\\\\
&=& D^*_N(\xi_1,\ldots,\xi_N).\cqfd
 \end{eqnarray*}

The above claim $(b)$ and the corollary both emphasize the fact that  considering uniform weights $\frac1N$ induces the loss of a $\log N$ factor compared to an optimal (or simply  rate optimal) greedy sequence   for optimal quantization since, for such an $(L^1,U([0,1]))$ greedy optimal sequence $a=(a_N)_{N\ge 1}$, one has for every integer $N\ge 1$,
\[
\left|\int_{[0,1]^d} f(u)du-  \sum_{i=1}^N w^{(N)}_i  f(a_i)\right|\le \kappa(a)\,[f]_{\rm Lip} \frac{1}{N^{\frac 1d}}
\]
where the $N$-tuple $(w^{(N)}_i)_{1\le i\le N}$ is vector of   hyper-volumes (Lebesgue measure) of the Voronoi cells attached to $a^{(N)}$.  Of course the practical implementation of such greedy sequences remains more demanding since one needs to have access to these $N$-tuples of weights. 

However, by contrast,    optimal quantization based  cubature formulas  
turn out to be  efficient (accurate) for much lower values of $N$ than sequences with low discrepancy (see $e.g.$ the numerical experiment carried out in~\cite{PAPR} dealing with the pricing of European derivatives).

\appendix
\section{Appendix: Greedy Lloyd's~I procedure}\label{app:A}
\subsection{The one-dimensional greedy Lloyd~I procedure} \label{app:A1}
The first is to establish the uniqueness of the equilibrium point $a_{_N}$ satisfying~\eqref{eq:Stationarity} and the convergence of the Lloyd~I procedure {\em at level $N$} toward this point, but  with the significant additional constraint that  the endpoints of the (closed convex) support of the strongly unimodal distribution $\mu$ are {\em active} (though fixed). By active we mean that, when finite, they have there own Voronoi area. To be more precise we will show the following proposition

\begin{Pro}\label{pro:lloydI-I} Let $\mu$ be a distribution on the real line  with $\log$-concave density $:varphi$ ($i.e.$ strongly unimodal) with a finite second moment. Then the support $I=\overline{\{\varphi>0\}})$   is  closed interval  with endpoints $a,\, b\!\in\overline{ \R}$. If $a$ or $b$ are finite, one may assume without loss of generality that $\varphi(a)$ or $\varphi(b)>0$ (so that $I$ is closed). Then the quantization problem at level $N$ with active  finite endpoints (if any) reads
\[
\min_{x\in I} \left[ \varphi(x) := \E \left(|X-a|\wedge|X-b|^2\wedge |X-x|^2\right)\right]
\] 
(note that when $a$ or $b$ are infinite, the corresponding terms in the above expectation can be omitted).  

\smallskip
\noindent $(a)$ The function $\varphi$ is differentiable on $I$ with a derivative given, for every $x\!\in I$, by 
\[
 \varphi'(x) = \frac 12\int_{\frac{a+x}{2}}^{\frac{x+b}{2}}(x-\xi)\mu(d\xi).
\]
Furthermore ${\rm argmin}_{I}\,{\cal G}$ is reduced to a single (stationary) point $x^*$ satisfying $\varphi'(x^*)=0$ $i.e.$
\[
x^*= \Phi(x^*)\quad \mbox{ where }\quad \Phi(x)= \frac{K_{\mu}\big(\frac{b+x}{2}\big)-K_{\mu}\big(\frac{a+x}{2}\big)}{F_{\mu}\big(\frac{b+x}{2}\big)-F_{\mu}\big(\frac{a+x}{2}\big)}
\]
and $F_{\mu}$ and $K_{\mu}$ denote the cumulative distribution and first moment functions of the distribution $\mu$ respectively.

\smallskip
\noindent $(b)$ The {\em greedy Lloyd~I} procedure  defined by 
\[
x_{n+1} = \Phi(x_n), \; x_0\!\in I
\]
converges toward $x^*$
\end{Pro}

This result can be seen as a variant of the Lloyd procedure at levels $N$ ($N=1$ up tp $3$), depending on the finiteness of the endpoints of the interval $I$.

\medskip
\noindent {\bf Proof.}  First note that,  when both endpoints are infinite and cannot be active, the above statement becomes trivial since ${\cal G}(x)= \E |X-x|^2$ which attains its minimum at $x^*= \E\,X$, whereas  the Lloyd~I procedure reads $x_1 = \E\,X$, $n\ge 1$, whatever the starting point $x_0$ is.

Otherwise, if $a$ or $b$ are finite, we may assume, up to a symmetry-translation, that   $a=0$ and $b\!\in (0,+\infty]$.
\smallskip
\noindent$(a)$ Elementary computations show that, for every $x\!\in I$, 
\[
{\cal G}(x) = \frac 12\left[ \int_0^{\frac x2}\xi^2\mu(d\xi) +\int_{\frac x2}^{\frac{x+b}{2}}(x-\xi)^2\mu(d\xi)+\int_{\frac{x+b}{2}}^{+\infty}(x-\xi)^2\mu(d\xi)\mbox{\bf 1}_{\{b<+\infty\}}  \right]
\]
\[
{\cal G}'(x)= \int_{\frac x2}^{\frac{x+b}{2}}(x-\xi)\mu(d\xi)
\]
and
\[
{\cal G}''(x)= F_{\mu}\big(\frac{x+b}{2}\big) -F_{\mu}\big(\frac{x}{2}\big) -\Big(\frac{x+b}{2}-x\Big) f\big(\frac{x+b}{2}\big)-\Big(x-\frac x2\Big)f\big(\frac x2\big).
\]
In what follows we focus on the case $b<+\infty$. The case $b=+\infty$ can be handled likewise (in fact in an easier way). 

Note that ${\cal G}'(0)= -\int_0^{\frac 12} \xi\mu(d\xi)<0$ and ${\cal G}'(b)= \int_{\frac 12}^1(b-\xi)\mu(d\xi)>0$ so that ${\cal G}$ has at least one zero on $(0,b)$. (When $b=+\infty$, the existence follows form the fact that ${\cal G}$ does attain a minimum on $(0, +\infty)$ since $\lim_{x\to +\infty}{\cal G}(x)=+\infty$.)

Set $y_1= \frac x2$ and $y_2= \frac{x+b}{2}$. If we assume that $x$ is a solution to $x=\Phi(x)$ (or equivalently to the stationary point equation ${\cal G}'(x)=0$), we can plug   this expression for $x$ into the above equation for ${\cal G}''(x)$ so that ${\cal G}''(x)$ can be expressed as a function of  $y_1$ and $y_2$ ias follows:
\[
{\cal G}''(x)=\frac{\widetilde \Phi(y_1,y_2)}{F_{\mu}(y_2) -F_{\mu}(y_1)},\; y_1<y_2, \;y_1,y_2\!\in I,
\]
with 
\[
\widetilde \Phi(y_1,y_2) = \big(F_{\mu}(y_2) -F_{\mu}(y_1)\big)^2+\big(K_{\mu}(y_2) - K_{\mu}(y_1)\Big)\big(\varphi(y_2)-\varphi(y_1)\big)- \big(F_{\mu}(y_2) -F_{\mu}(y_1)\big)\big(y_2\,\varphi(y_2)-y_1\,\varphi(y_1)\big).
\]
Now we   consider $y_1$ and $y_2$ as free variables living in $I$ such that $y_1\le y_2$. First we note that $\widetilde \Phi(y,y)=0$. Then, denoting by $\varphi'_r$ the right derivative of the $\log$-concave function $\varphi$, we compute the following two (right) partial derivatives of $\widetilde \Phi$: 
\begin{eqnarray*}
\left(\frac{\partial \widetilde \Phi}{\partial y_1}\right)_{\!\!r}\hskip -0.15 cm(y_1,y_2)&=&\big(F_{\mu}(y_2) -F_{\mu}(y_1)\big)\big(y_1\varphi'_r(y_1)-\varphi(y_1)\big)+(y_2-y_1)\varphi(y_1)\varphi(y_2)-\varphi'_r(y_1)\big(K_{\mu}(y_2) - K_{\mu}(y_1)\big)\quad \\
\mbox{and }\hskip 2 cm  &&\\ 
\left(\frac{\partial^2 \widetilde \Phi}{\partial y_1\partial y_2}\right)_{\!\!r}\hskip -0.15 cm (y_1,y_2)&=& (y_2-y_1) \Big(\varphi(y_1)\varphi'_r(y_2)-\varphi'_r(y_1)\varphi(y_2)\Big).
\end{eqnarray*}
As $(\log \varphi)'_r = \frac{\varphi'_r}{\varphi}$ is non-increasing in $I$, it follows that $\frac{\partial^2 \widetilde \Phi}{\partial y_1\partial y_2}(y_1,y_2)<0$ if $y_1<y_2$ so that 
$y_2\mapsto \Big(\frac{\partial \widetilde \Phi}{\partial y_1}\Big)_r(y_1,y_2)$  is (strictly) decreasing on $[y_1,b)$ which in turn implies it is positive on $(y_1, b)$. This shows that $\Phi(y_1,y_2)>0$ for every $y_1,y_2\!\in I$, $y_1<y_2$. As a consequence, any stationary point $x$ satisfies ${\cal G}''(x)<0$ $i.e.$ is a strict local minimum of ${\cal G}$. This implies uniqueness of the solution to the equation ${\cal G}'(x)=0$ by an elementary  one dimensional ``mountain pass" argument.

\smallskip
\noindent $(b)$ In this second claim,  we use again  a random variable $X$ with distribution $\mu$. By Proposition~\ref{pro:greedex}$(b)$, we know that if $W_{[n]}= [\frac{x_n}{2}, \frac{x_n+b}{2}]$ denotes the closed Voronoi cell of $x_n$ with respect to $\{0,x_n,b\}$ (if $b$ is finite, or $\{0,x_n\}$ otherwise) then 
\[
\E_{\mu} \big(|X-x_{n+1}|^2\mbox{\bf 1}_{\{X\in W_{[n]}\}}\big)\le  \E_{\mu} \big(|X-x_{n}|^2\mbox{\bf 1}_{\{X\in W_{[n]}\}}\big).
\]
with equality iff $x_{n+1}= x_n$ which is equivalent to $x_n= x^*$ (see claim~$(a)$ above). Decomposing $|X-x_{n+1}|^2$ on the Voronoi partition $[0, \frac{x_n}{2}]\cup[\frac{x_n}{2}, \frac{x_n+b}{2}]\cup[\frac{x_n+b}{2},b]$ of $I$, one derives that ${\cal G}(x_{n+1})<{\cal G}(x_n)$ as soon as $x_n\neq x^*$.  the function ${\cal G}$ being non-negative ${\cal G}(x_n)\to \ell$ as $n\to \infty$.

\smallskip When $b$ is finite the sequence $(x_n)_{n\ge 0}$ is trivially bounded. When $b=+\infty$, assume there exists a subsequence $x_{n'}\to +\infty$. By combining the above monotony property and Fatou's Lemma, we get 
\[
\E|X|^2\wedge|X-x_0|^2= {\cal G}(x_0)\ge \liminf_n {\cal G}(x_{n'})\ge \E\, X^2
\]
which implies that $X^2\le (X-x_0)^2$ $\P$-$a.s.$ This  is clearly not satisfied on the event  $\{X\!\in [0, \frac{x_0}{2})\}$ which has positive  probability. Consequently, $(x_n)_n$ is always bounded.

 Then let $x_{\infty}=\lim_{n\to +\infty} x_{n'}$   be a limiting value of the $I$-valued sequence $(x_n)_{n\ge 0}$. Up to a new extraction, still denoted $(n')$, one may  assume that $x_{n'+1}$ converges toward a limiting value $x'_{\infty}$ as well. Passing to the limit owing to continuity we get 
 $$
 x'_{\infty} = \int_{\frac{x_{\infty}}{2}}^{\frac{x_{\infty}+b}{2}}\xi\mu(d\xi).
 $$
One shows as above that,  except if $x'_{\infty}=x_{\infty}=x^*$, ${\cal G}(x'_{\infty})<{\cal G}(x_{\infty})$ which cannot be true since the sequence $({\cal G}(x_n))_{n\ge 0}$ converges to $\ell^2$. Consequently, $x^*=x_{\infty}$ is the only possible limiting value for the bounded sequence $(x_n)_{n\ge 0}$ $i.e.$ its limit.~$\cqfd$

\bigskip
\noindent {\bf Proof of Proposition~\ref{pro:Greedy-Lloyd}.} The result follows by applying the above result to the procedure on the interval $\displaystyle \left[\frac{a^{(N-1)}_{i_0}}{2}, \frac{a^{(N-1)}_{i_0+1}}{2}\right]$ of maximal inertia. $\cqfd$

\bigskip
 \noindent  {\bf Remark.}  If we  choose $a_{[0]}$ inside an interval which {\em has not}  the highest local inertia, the procedure will still converge since we never use this fact throughout the proof of the convergence. The resulting limit will live  in the same interval as the starting  value since the algorithm  leaves each interval stable by an obvious convexity argument. So the greedy Lloyd~I procedure yields potentially $N+1$ ``candidates" corresponding to each possible starting interval, but only one (issued from the interval with the highest local inertia) is solution to the greedy optimal quantization problem.


%
%

\subsection{The multi-dimensional greedy Lloyd~I procedure (proof of Proposition~\ref{pro:LloydI-d})}\label{app:A2} We assume in this section  that $\mu$ has a convex support $C_{\mu}={\rm supp}(\mu)$ and that $d\ge 2$. Note that in such  a framework there is a major topological difference with the $1$-dimensional case: a convex set not reduced to a single point remains pathwise connected when one point of its points is removed. Owing to that property, it is easy to show that  the algorithm may visit with positive probability the whole support of $C_{\mu}$ (to be precise any nonempty open set of $C_{\mu}$). Moreover,  the points can no longer be naturally ordered like  in $1$-dimension.

To alleviate notations, we  denote by ${\cal G}$ the $\R_+$-valued function $a\mapsto e_{2}\big(a^{(N-1)}\cup\{a\}\big)^2 $ defined on $C_{\mu}$ by 
\[
{\cal G}(a) = e_{2}\big(a^{(N-1)}\cup\{a\}\big)^2= \E \Big(d\big(X, a^{(N-1)}\cup \{a\}\big)^2\Big).
\]
Let $a_{[0]}\!\in C_{\mu} \setminus a^{(N-1)}$. Lloyd's~I procedure is defined by induction by Equation~\eqref{eq:Lloyd-I-d}, namely 
\[
a_{[n+1]} = \E \big(X\,|\, X\!\in W_{N,[n]}\big)\!\in C_{\mu}
\] 
where  $W_{N,[n]}$ denotes the (closed) Voronoi cell of $a_{[n]}$ induced by $a^{(N-1)}\cup\{a_{[n]}\}$. 
 
\noindent {\sc Step~1:} It follows from Proposition~\ref{pro:greedex}$(b)$ that, as son as $a_{[n]}$ is not stationary, $i.e.$ $a_{[n]}\neq \E \big(X\,|\, X\!\in W_{N,[n]}\big)$, one has 
 \[
 \E\big( |X-a_{[n+1]}|^2\mbox{\bf 1}_{\{ X \in W_{N,[n]} \}} \big)<  \E\big( |X-a_{[n]}|^2 \mbox{\bf 1}_{\{ X \in W_{N,[n]} \} }\big).
 \]
 hence
 \begin{eqnarray*}
 {\cal G}(a_{[n+1]})&=& \E \big(d(X,a^{(N-1)}\cup\{a_{[n+1]}\})^2\big) \\
 &\le&   \E\big( d(X,a^{(N-1)})^2\mbox{\bf 1}_{\{ X\notin W_{N,[n]}\}}\big)+ \E\big( |X-a_{[n+1]}|^2\mbox{\bf 1}_{\{ X\in W_{N,[n]}\}}\big)\\
 &<&  \E\big( d(X,a^{(N-1)})^2\mbox{\bf 1}_{\{ X\notin W_{N,[n]}\}}\big)+ \E\big( |X-a_{[n]}|^2\mbox{\bf 1}_{\{ X\in W_{N,[n]}\}}\big)\\
 &=&{\cal G}(a_{[n]}).
 \end{eqnarray*}
  hence, the (non-increasing,  non-negative) sequence $({\cal G}(a_{[n]}))_n$ converges to a finite limit $\ell^2\!\in \R_+$ as $n\to +\infty$. The fact that $\ell\!\in \big(e_2(a^{(N)}, e_2(a^{(N-1)}\cup\{a_{[0]}\}\big)$ is obvious form what precedes.

 \noindent {\sc Step~2:} Assume there exists a subsequence $(a_{[n']})$ such that $|a_{[n']}|\to +\infty$ as $n\to +\infty$.  Combining the above monotony of the sequence $({\cal G}(a_{[n]}))_{n\ge 0}$ and Fatou's Lemma yields
 \[
 \E \,d(X,a^{(N-1)})^2 \le \liminf_n \E \,d(X,a^{(N-1)}\cup\{\a_{[n']}\})^2\le \liminf_n {\cal G}(a_{[n]})\le {\cal G}(a_{[0]})
 \]
But, as $a_{[0]}\!\in C_{\mu}\setminus a^{(N-1)}$, we know from Proposition~\ref{pro:greedex}$(a)$  that  $ {\cal G}(a_{[0]})<\E \,d(X,a^{(N-1)})^2 $ which yields a contradiction.
 
 \medskip
 \noindent {\sc Step~3}: Let $a_{[\infty]}$ be a limiting value of the bounded sequence $(a_{[n]})_n$ ($i.e.$ the limit of a subsequence). Up to a new extraction, we may also assume that $a_{[n+1]}\to a'_{[\infty]}$. Since ${\cal G}(a_{[n]})$ is non-decreasing, ${\cal G}( a_{[\infty]})$ and ${\cal G}(a'_{[\infty]})\le{\cal G}(a_{[0]})<e_{2}(a^{(N-1)},\mu)$ so that we $ a'_{[\infty]}$, $ a'_{[\infty]}\notin a^{(N-1)}$. The distribution $\mu$ being strongly continuous ($i.e.$ assigning no mass to hyperplanes), one shows $e.g.$ by  following the lines of the proof of Lemma~2.3  in~\cite{PAYU}  that 
\[
\E\big(X\,|\, W\!\in W_{N,[n]}\big)\to \E \big(X\,|\, X\!\in W_{N,[\infty]}\big)
\]
where $W_{N,[\infty]}$ denotes the closed Voronoi cell of $a_{[\infty]}$ induced by $a^{(N-1)}\cup\{a_{[\infty]}\}$. Consequently
\[
a'_{\infty} =  \E \big(X\,|\, X\!\in W_{N,[\infty]}\big). 
\]
If $a_{[\infty]}\neq a_{[\infty]}$ then ${\cal G}(a'_{[\infty]})< {\cal G}(a_{[\infty]})$ which is in a contradiction with the fact  $({\cal G}(a_{[n]}))_n$ converges to a finite limit $\ell^2\!\in \R_+$ as $n\to +\infty$. Hence $a'_{[\infty]}=a_{[\infty]}$ which shows that one the one hand that $a_{[n+1]}-a_{[n]}\to 0$ as $n\to +\infty$ and that any limiting value of $(a_{[n]})_n$ is a stationary point in the sense that $\displaystyle a_{\infty} =  \E \big(X\,|\, X\!\in W_{N,[\infty]}\big)$.  The conclusion follows by standard topological arguments on convergence of sequences.$\cqfd$

\section{Appendix: A technical result on sequences}\label{app:B}
\begin{Lem} \label{lem:A-N} Let $(A_N)_{N\ge 1}$ be a   sequence of non-negative real numbers   and let $\rho\!\in (0,+\infty)$ such that 
\[
\forall\, N\ge 1, \; A _{N+1} \le A _N-C\,A_N^{1+\rho}
\]
for some real constant $C>0$. Then there exists a real constant $K>0$ such that
\[
\forall\, N\ge 1,\quad A_N\le  K N^{-\frac{1}{\rho}}.
\]
\end{Lem}

\noindent {\bf Proof.} We may assume that $A_{_N}>0$ for every $N\ge1$, it follows from the inequality satisfies by the sequence $(A_N)_{N\ge 1}$ that for very $N\ge 1$, 
\[
\frac{1}{A^{\rho}_{N+1}}\ge \frac{1}{A^{\rho}_N}\frac{1}{(1-CA^{\rho}_N)^{\rho}}\ge \frac{1}{A^{\rho}_N}(1+CA^{\rho}_N)^{\rho}
\]
Now, there exists $u_0=u_0(\rho)$ such that for every $u\!\in [0,u_0]$, 
$(1+u)^{\rho }\ge 1+ \frac{\rho}{2}u$. It is clear from the assumptions that $A_N\downarrow 0$, hence, there exists  a large enough  integer $N_0$  such that for every $N\ge N_0$, 
\[
\frac{1}{A^{\rho}_{N+1}}\ge \frac{1}{A^{\rho}_N} + \frac{C\rho}{2}
\]
which in turn implies that
\[
\frac{1}{A^\rho_{N}}\ge \frac{1}{A^\rho_{N_0}} + \frac{C\rho}{2}(N-N_0)\ge  \frac{C\rho}{2}(N-N_0)
\] 
so that, for every $N>N_0$, 
\[
A_N \le \left(\frac{2}{C\rho}\right)^{\frac{1}{\rho}} \frac{1}{(N-N_0)^{\frac{1}{\rho}}}.
\]
This  completes the proof.~$\cqfd$

\end{document}